\documentclass[draft,11pt,english]{article}

\usepackage{amsfonts,amsmath,amsthm,amscd,amssymb,latexsym,verbatim,
caption2}

\theoremstyle{plain}
\newtheorem{theorem}{Theorem}[section]

\newtheorem{proposition}{Proposition}[section]
\newtheorem{corollary}{Corollary}[section]
\newtheorem{definition}{Definition}[section]
\theoremstyle{definition}

\newtheorem{remark}{Remark}[section]
\newcommand{\keywords}{\textbf{Key words. }\medskip}
\newcommand{\subjclass}{\textbf{MSC 2010. }\medskip}
\renewcommand{\abstract}{\textbf{Abstract. }\medskip}

\numberwithin{equation}{section}

\begin{document}

\title{\textbf{Rational Matrix-Valued Pick Functions of Several Variables}}

\author{M.F. Bessmertny\u{\i}}

\date{}\maketitle

\begin{abstract}
For a rational function of several variables with nonnegative imaginary part on the upper poly-half-plane, the matrix representations are obtained.
\end{abstract}
\medskip



\section{Introduction}\label{s1:1}

  Let $\Pi^{d}=\{z\in\mathbb{C}^{d}\mid \mbox{Im}\,z_{1}>0,\ldots,\mbox{Im}\,z_{d}>0\}$ be an open upper poly-half-plane.
    \emph{The Pick class}
    $\mathcal{P}_{d}^{m\times m}$ consists of $\mathbb{C}^{m\times m}$-valued functions $f(z_{1}.\ldots,z_{d})$ holomorphic in $\Pi^{d}$ and satisfying the condition
      \begin{equation}\label{eq1.1}
        (f(z)-f(z)^{\ast})/2i\geq0\quad \mbox{for} \quad z\in\Pi^{d},
      \end{equation}
  where $^{\ast}$ means transition to the Hermitian conjugate matrix.
  We will consider only \emph{rational} function of the Pick class, as well as a subclass $\mathcal{IP}_{d}^{m\times m}$ of rational Pick functions, which are \emph{Cayley inner}.
  The latter means that such a function takes Hermitian matrix values almost everywhere on distinguished boundary $\mathbb{R}^{d}$ of the upper poly-half-plane $\Pi^{d}$,
  or equivalently, is the double Cayley transform (over the variables and over the matrix values) of an inner
  function on the unit polydisk  $\mathbb{D}^{d}=\{z\in\mathbb{C}^{d}\mid |z_{1}|<1,\ldots,|z_{d}|<1\}$.
  The structure of scalar rational inner functions in polydisk is considered in \cite{uj22}.

  For one-variable functions $f(\lambda)$ of the Pick class\footnote{Other popular titles for the same class of
  one-variable functions are ``Nevanlinna", ``Nevanlinna-Pick", and ``$R$-functions".},
  Nevanlinna's integral representation is well know. Pick's functions of several variables have been studied by many authors (see, for example, bibliography in \cite{uj01}).
  In \cite{uj01} for the L\"{o}wner class a generalisation of the classical Nevanlinna theorem to several variables was considered.
  The L\"{o}wner class $\mathcal{L}_{d}$ is a set of scalar Pick functions $h$ in $d$ variables such that
      \begin{equation}\label{eq1.2}
      h(z)-\overline{h(w)}=\sum_{k=1}^{d}(z_{k}-\overline{w}_{k})\Phi_{k}(z,w),
      \end{equation}
 where $\Phi_{k}(z,w)$ is positive semidefinite, that is, for all $n\geq1$, $z_{1},\ldots,z_{n}\in\Pi^{1}$, $c_{1},\ldots,c_{n}\in\mathbb{C}$,
      $$
      \sum_{i,j=1}^{n}\Phi_{k}(z_{j},z_{i})\overline{c}_{i}c_{j}\geq0.
      $$
 Note the checking of condition (\ref{eq1.2}) for a function is a difficult.

 The Herglots-Agler class consists of the functions which are holomorphic on the right poly-half-plane and whose values on any commutative $d$-tuple
 of strictly accretive operators on Hilbert space have positive semidefinite real part.
 For rational Cayley inner functions $f(z)$ of the Gerglots-Agler class,
 \emph{a long-resolvent representation} was obtained in \cite{uj02}:
       \begin{equation}\label{eq1.3}
       f(z)=A_{11}(z)-A_{12}(z)A_{22}(z)^{-1}A_{21}(z),
       \end{equation}
 where
       $$
         \begin{pmatrix}
         A_{11}(z) & A_{12}(z) \\
         A_{21}(z) & A_{22}(z)
         \end{pmatrix}=
            H+z_{1}A_{1}+\cdots+z_{d}A_{d},
       $$
  and matrix coefficients satisfy the conditions
          \begin{equation}\label{eq1.4}
          H^{\ast}=-H,\quad A_{k}\geq0,\,k=1,\ldots,d.
          \end{equation}
  It is easy to see that if in (\ref{eq1.4}) the matrix $H$ satisfies condition
       \begin{equation}\label{eq1.5}
       (H-H^{\ast})/2i\geq0,
       \end{equation}
  then (\ref{eq1.3}) is a rational matrix-valued function of the Pick class. Indeed, since $A_{22}(z)$ is an invertible matrix, we see that for $z,\zeta\in\mathbb{C}^{d}$:
      $$
      f(z)=
      \begin{pmatrix}
      I_{m} \\
      -A_{22}(\zeta)^{-1}A_{21}(\zeta)
      \end{pmatrix}^{\ast}
          \begin{pmatrix}
          A_{11}(z) & A_{12}(z)\\
          A_{21}(z) & A_{22}(z)
          \end{pmatrix}
               \begin{pmatrix}
               I_{m} \\
               -A_{22}(z)^{-1}A_{21}(z)
               \end{pmatrix}.
      $$
  From here we get
         \begin{multline}\label{eq1.6}
         \qquad\frac{f(z)-f(\zeta)^{\ast}}{2i}=\\
         (I_{m},\,-A_{21}(\zeta)^{\ast}A_{22}(\zeta)^{-1\ast})\frac{H-H^{\ast}}{2i}
           \begin{pmatrix}
                         I_{m}  \\
             -A_{22}(z)^{-1}A_{21}(z)
           \end {pmatrix}+\\
               \sum_{k=1}^{d}\frac{z_{k}-\overline{\zeta}_{k}}{2i}(I_{m},\,-A_{21}(\zeta)^{\ast}A_{22}(\zeta)^{-1\ast})A_{k}
                    \begin{pmatrix}
                         I_{m}  \\
                    -A_{22}(z)^{-1}A_{21}(z)
                    \end{pmatrix}.
         \end{multline}
  For $\zeta=z$ from (\ref{eq1.6}) we obtain $f(z)\in\mathcal{P}_{d}^{m\times m}$.
  \begin{remark}\label{rem1.1} Let be $(H-H^{\ast})/2i=B_{0}^{\ast}B_{0}$, $A_{k}=B_{k}^{\ast}B_{k}$, $k=1,\ldots,d$.
          Consider rational functions
        $$
        \theta_{k}(z)=B_{k}
            \begin{pmatrix}
            I_{m}  \\
            -A_{22}(z)^{-1}A_{21}(z)
            \end{pmatrix},\quad k=0,1,\ldots,d.
        $$
      Then from (\ref{eq1.6}) we obtain the positive-kernel decomposition
        $$
        \frac{f(z)-f(\zeta)^{\ast}}{2i}=\theta_{0}(\zeta)^{\ast}\theta_{0}(z)+\sum_{k=1}^{d}\frac{z_{k}-\overline{\zeta}_{k}}{2i}\theta_{k}(\zeta)^{\ast}\theta_{k}(z).
        $$
   \end{remark}
  \medskip

  Another way to obtain rational functions of the Pick class is as follows. Let
          $$
          H=
          \begin{pmatrix}
          A & B \\
          C & D
          \end{pmatrix},\quad \frac{H-H^{\ast}}{2i}\geq 0
          $$
   be a $z$-independent block matrix, where block $D$ has size $n\times n$.
   If
       $$
       Z_{n}=\text{diag}\,\{I_{n_{0}},z_{1}I_{n_{1}},\ldots,z_{d}I_{n_{d}}\},\quad n_{0}+n_{1}+\cdots+n_{d}=n,
       $$
   then function
          \begin{equation}\label{eq1.7}
          f(x)=A+BZ_{n}(I_{n}-DZ_{n})^{-1}C
          \end{equation}
   is the rational matrix-valued function of the Pick class.
   Indeed,  $f(z)=A+BZ_{n}(I-DZ_{n})^{-1}C=A-B(D-Z_{n}^{-1})^{-1}C$. Then
             $$
             f(z)=
                 \begin{pmatrix}
                        I     \\
                 -(D-Z_{n}^{-1})^{-1}C
                 \end{pmatrix}^{\ast}
                     \begin{pmatrix}
                      A   &    B \\
                      C   & D-Z_{n}^{-1}
                     \end{pmatrix}
                          \begin{pmatrix}
                               I     \\
                          -(D-Z_{n}^{-1})^{-1}C
                          \end{pmatrix}.
             $$
  From here
        \begin{multline}\label{eq1.8}
          \frac{f(z)-f(z)^{\ast}}{2i}=
                 \begin{pmatrix}
                        I     \\
                 -(D-Z_{n}^{-1})^{-1}C
                 \end{pmatrix}^{\ast}
                      \frac{H-H^{\ast}}{2i}
                          \begin{pmatrix}
                               I     \\
                          -(D-Z_{n}^{-1})^{-1}C
                          \end{pmatrix}+\\
        C^{\ast}(I-DZ_{n})^{-1\ast}\frac{Z_{n}-Z_{n}^{\ast}}{2i}(I-DZ_{n})^{-1}C.\qquad
        \end{multline}
   It follows from (\ref{eq1.8}) that $f(z)\in\mathcal{P}_{d}^{m\times m}$.

   For functions of the Pick class, the representation (\ref{eq1.7}) is analog of the transfer-function realisation.
    \medskip

   In this article, we will show that functions of the form (\ref{eq1.7}) and (\ref{eq1.3}) (subject to condition (\ref{eq1.5}))
   exhaust rational $\mathbb{C}^{m\times m}$-valued functions of the Pick class $\mathcal{P}_{d}^{m\times m}$:

\begin{theorem}\label{th7.2} \textbf{\emph{(Main Theorem)}}.
       Let $f(z)$ be a $\mathbb{C}^{m\times m}$valued function of $d$ complex variables. The following statements are equivalent.
     \medskip

    \textbf{\emph{(0)}} $f(z)$ is a rational Pick function.
     \medskip

    \textbf{\emph{(1)}} There exist $n, n_{1},\ldots,n_{d}\in\mathbb{Z}_{+}$, $ n_{0}\geq0$, $n=n_{0}+n_{1}+\cdots+n_{d}$
         and the $z$-independent block $(m+n)\times (m+n)$ matrix
             \begin{equation}\label{eq7.21}
                H=
                 \begin{pmatrix}
                   A  &  B  \\
                   C  &  D
                 \end{pmatrix},\quad
                      \frac{H-H^{\ast}}{2i}\geq0,
             \end{equation}
     with block $D$ of size $n\times n$ such that
              \begin{equation}\label{eq7.22}
                f(z)=A-B(D+Z_{n})^{-1}C,
              \end{equation}
    where $Z_{n}=\emph{diag}\{0_{n_{0}},z_{1}I_{n_{1}},\ldots,z_{d}I_{n_{d}}\}$.
    \medskip

    \textbf{\emph{(2)}} There exist $n, n_{1},\ldots,n_{d}\in\mathbb{Z}_{+}$, $ n_{0}\geq0$, $n=n_{0}+n_{1}+\cdots+n_{d}$
        and the $z$-independent block  $(m+n)\times (m+n)$ matrix
             \begin{equation}\label{eq7.22}
                H=
                 \begin{pmatrix}
                   A  &  B  \\
                   C  &  D  \\
                 \end{pmatrix},\quad
                      \frac{H-H^{\ast}}{2i}\geq0,
             \end{equation}
    with block $D$ of size $n\times n$ such that
              \begin{equation}\label{eq7.23}
                f(z)=A+BZ_{n}(I_{n}-DZ_{n})^{-1}C,
              \end{equation}
    where $Z_{n}=\emph{diag}\{I_{n_{0}},z_{1}I_{n_{1}},\ldots,z_{d}I_{n_{d}}\}$.
    \medskip

    \textbf{\emph{(3)}} The function $f(z)$ has a long-resolvent representation
             \begin{equation}\label{eq7.24}
             f(z)=A_{11}(z)-A_{12}(z)A_{22}(z)^{-1}A_{21}(z),
             \end{equation}
    where
             \begin{equation}\label{eq7.25}
                 \begin{pmatrix}
                 A_{11}(z) & A_{12}(z) \\
                 A_{21}(z) & A_{22}(z)
                 \end{pmatrix}=
                     H+z_{1}A_{1}+\cdots+z_{d}A_{d},
             \end{equation}
    and the matrix coefficients satisfy the conditions:
               \begin{equation}\label{eq7.26}
               (H-H^{\ast})/2i\geq0,\quad A_{k}\geq0 \quad\text{for}\quad k=1,\ldots,d.
               \end{equation}

  \noindent If $f(z)$ is an inner function, then the matrices $H$ in \textbf{\emph{(1), (2), (3)}} can be chosen as Hermitian \emph{(}$H=H^{\ast}$\emph{)}.
       \end{theorem}

 The article is organized as follows. In Section \ref{s:2}, we reduce the study of rational Pick function
 to the study of a multi-affine function of the Pick class, the degree of which in each variable does not exceed 1.

 Section \ref{s:3} proved (Theorem \ref{th3.1}) that every rational Pick function in $d$ variables can be obtained
 from the Cayley inner function in $d+1$ variables of the Pick class. In Theorem \ref{th3.2} (based on Theorem \ref{th3.3}), a criterion is obtained
 for a multi-affine function to be an inner function of the Pick class.

 The necessary apparatus for obtaining representations is developed in Sections \ref{s:4}, \ref{s:5} and \ref{s:6}.

 The basis is the Sum-of-Squares Theorem (Theorem \ref{th4.1}).
 The proof of the Sum-of Squares Theorem is given in  Appendix (Section \ref{s:Ap}).

 In Section \ref{s:5}, for multi-affine functions of the Pick class in several variables,
 a generelization of the Darlington method of realization a function is obtained (Theorem \ref{eq5.1}).

 A superposition of matrix coefficients of fractional linear transformations of a special form is considered in Proposition \ref{pr6.1}.

 The matrix representation of Cayley inner functions of the Pick class was obtained in Section \ref{s:7} (Theorem \ref{th7.1}).
 The proof of the Main Theorem is given in Section \ref{s:8}.

 Preliminary information is provided in each section as needed.

  \section{Reduction to the Multi-Affine Functions}\label{s:2}

  The rational $\mathbb{C}^{m\times m}$-valued function will be written in the form $f(z)=P(z)/q(z)$,
  where $P(z)$ is a $\mathbb{C}^{m\times m}$-valued polynomial and
  $q(z)$ is a scalar $\mathbb{C}$-valued polynomial. In fact, division $P(z)/q(z)$ is the standard operation of multiplying of the matrix $P(z)$  by the number $q(z)^{-1}$.

   We say that the polynomial $p(z)$ \emph{affine} in $z_{k}$ if $\mbox{deg}_{z_{k}}p(z)=1$, and we say $p(z)$ is \emph{multi-affine},
   if it is affine in $z_{k}$ for all $k=1,\ldots,d$. A rational function $P(z)/q(z)$ will be called \emph{multi-affine} if
       $$
       \max\{\mbox{deg}_{z_{k}}P(z),\,\mbox{deg}_{z_{k}}q(z)\}=1,\quad \forall k=1,\ldots,d.
       $$
  Elementary symmetric polynomial are examples of multi-affine functions:
             $$
              \sigma_{k}(\zeta_{1},\ldots,\zeta_{n})=\sum_{i_{1}<i_{2}<\cdots <i_{k}}\zeta_{i_{1}}\zeta_{i_{2}}\cdots \zeta_{i_{k}},\, 0<k\leq n,
              \quad \sigma_{0}(\zeta_{1},\ldots,\zeta_{n})\equiv 1.
             $$

    Reduction of the problem to the multi-affine case is based on the use of the degree reduction operator \cite{uj12,uj18}.

   \begin{definition}\label{def2.1}
       \emph {Let $p(z_{0},z)=\sum_{k=0}^{n_{0}}p_{k}(z) z_{0}^{k}$ be a polynomial and $n_{0}\leq n\in\mathbb{N}$. A map}
           \begin{equation}\label{eq2.1}
             \mathbf{D}_{z_{0}}^{n}: \sum_{k=0}^{n_{0}}p_{k}(z)z_{0}^{k}\,\mapsto \,\sum_{k=0}^{n_{0}}p_{k}(z)\binom {n}k ^{-1}
             \sigma_{k}(\zeta_{1},\ldots,\zeta_{n})
           \end{equation}
       \emph {is called} a degree reduction operator \emph{in the variable $z_{0}$. If
      $f(z_{0},z)=P(z_{0},z)/q(z_{0},z)$ is rational function and $\mbox{deg}_{z_{0}}f(z_{0},z)=n_{0}$, then the degree reduction operator is defined as}
                \begin{equation}\label{eq2.2}
                  \mathbf{D}_{z_{0}}^{n}[\,P(z_{0},z)/q(z_{0},z)\,]:=\mathbf{D}_{z_{0}}^{n}[\,P(z_{0},z)\,]/\mathbf{D}_{z_{0}}^{n}[\,q(z_{0},z)\,].
                \end{equation}
   \end{definition}

  Under the condition $\zeta_{1}=\cdots =\zeta_{n}=z_{0}$ we get the original function. Thus, the operator $\mathbf{D}_{z_{0}}^{n}$ is invertible.
  It turns out that the degree reduction operator (\ref {eq2.2}) has the following important property.
    \begin{theorem}\label{th2.1}
      Let $f(z_{0},z)=P(z_{0},z)/q(z_{0},z)\in\mathcal{P}_{d+1}^{m\times m}$ and $P,q$ are coprime polynomials.
       If $\emph{deg}_{z_{0}}f(z_{0},z)=n_{0}$ and $n_{0}\leq n\in\mathbb{N}$, then
      $\widehat{f}(\zeta_{1},\ldots,\zeta_{n},z)=\mathbf{D}_{z_{0}}^{n}[\,f(z_{0},z)\,]$ is the Pick function of the class
      $\mathcal{P}_{d+n}^{m\times m}$, affine and symmetric in variables $\zeta_{1},\ldots,\zeta_{n}$, Moreover,
            \begin{equation}\label{eq2.3}
            \widehat{f}(z_{0},\ldots,z_{0},z)=f(z_{0},z).
            \end{equation}
    \end{theorem}
  We need the following statement.
  \medskip

  \noindent \textbf{Theorem (Grace-Walch-Szeg\"{o}).} (see \cite{uj12}, Theorem 2.12). \emph{Let $p$ be a symmetric multi-affine polynomial in $n$
  complex variables, let $\mathcal{C}$ be an open or closed circular region in $\mathbb{C}$, and let $\zeta_{1},\ldots,\zeta_{n}$ be any
  fixed points in the region $\mathcal{C}$, If $\emph{deg}\,p=n$ or $\mathcal{C}$ is convex, then there exists at least one point $\xi\in\mathcal{C}$
  such that $p(\zeta_{1},\ldots,\zeta_{n})=p(\xi,\ldots,\xi)$.}

  \begin{remark}
  Recall that a \emph{circular region} is a proper subset of the complex plane, which is bounded by circles (straight lines).
   In particular, the half-plane is a convex circular region.
  \end{remark}

    \emph{Proof of Theorem} \ref{th2.1}.
      A matrix-valued function $f(z_{0},z)$ is the Pick function if and only if for any row vector $\eta\in\mathbb{C}^{m}$ scalar function $\eta f(z_{0},z)\eta^{\ast}$
      is Pick's function. The coprime numerator and denominator of the scalar Pick function $p(z_{0},z)/q(z_{0},z)$ do not vanish $\Pi^{d+1}$.
      Since addition does not deduce from the class of Pick function, we see that $z_{d+1}+p(z_{0},z)/q(z_{0},z)$ is
      the Pick function in variables $z_{0},z,z_{d+1}$. Then its numerator satisfies the condition
        \begin{equation}\label{eq2.4}
          z_{d+1}q(z_{0},z)+p(z_{0},z)\neq0 \quad \mbox{for}\quad \mbox{Im}\,z_{0}>0,\, z\in\Pi^{d},\, \mbox{Im}\,z_{d+1}>0.
        \end{equation}
      Let us proove that the affine and symmetric in the variables $\zeta_{1},\ldots,\zeta_{n}$ polynomial
           \begin{multline}\label{eq2.5}
             \widetilde{p}(\zeta_{1},\ldots,\zeta_{n},z,z_{d+1})=\mathbf{D}_{z_{0}}^{n}[z_{d+1}q(z_{0},z)+p(z_{0},z)]=\\
            z_{d+1}\mathbf{D}_{z_{0}}^{n}[q(z_{0},z)]+\mathbf{D}_{z_{0}}^{n}[p(z_{0},z)],\quad
           \end{multline}
      do not vanish at
             \begin{equation}\label{eq2.6}
               \mbox{Im}\,\zeta_{j}>0,\, j=1,\ldots,n, \quad z\in\Pi^{d}, \quad \mbox{Im}\,z_{d+1}>0.
             \end{equation}
      Indeed, if the variables satisfy (\ref{eq2.6}) and $z, z_{d+1}$ is fixed, then by the Grace-Walsh-Szeg\"{o} Theorem there is $\xi,\,\mbox{Im}\,\xi>0$ such that
           $$
            \widetilde{p}(\zeta_{1},\ldots,\zeta_{n},z,z_{d+1})=\widetilde{p}(\xi,\ldots,\xi,z,z_{d+1})=z_{d+1}q(\xi,z)+p(\xi,z)\neq0.
           $$
      Let $\mbox{Im}\,\zeta_{j}>0,\, j=1,\ldots,n$, $z\in\Pi^{d}$ be fixed. Since the polynomial (\ref{eq2.5}) vanishes at the point
      $z_{d+1}=-\mathbf{D}_{z_{0}}^{n}[p(z_{0},z)]/\mathbf{D}_{z_{0}}^{n}[q(z_{0},z)]$ from the closed lower half-plane, we see that
              $$
                \mbox{Im}\,\widehat{f}(\zeta_{1},\ldots,\zeta_{n},z)=\mbox{Im}\,
                \left(\mathbf{D}_{z_{0}}^{n}[p(z_{0},z)]/\mathbf{D}_{z_{0}}^{n}[q(z_{0},z)]\right)\geq0.
              $$
      The relation (\ref{eq2.3}) is obvious
   \qed

 \section{Cayley inner functions of the Pick class}\label{s:3}

  We will reduce the problem of obtaining a representation of an arbitrary rational Pick function to the problem
  of representing a rational Pick function, which is Cayley inner.
  Then, by Theorem \ref{th2.1}, it suffices to restrict ourselves to multi-affine Cayley inner Pick functions.

  In addition, we obtain a criterion (Theorem \ref{th3.2}) for a multi-affine function to belong to the subclass
  $\mathcal{IP}_{d}^{m\times m}$ of Cayley inner functions of Pick class $\mathcal{P}_{d}^{m\times m}$.

    At the points of continuity, rational Cayley inner function $f(z)=P(z)/q(z)$ of the Pick class satisfies the condition
        \begin{equation}\label{eq3.1}
        f(\overline{z})^{\ast}= f(z).
        \end{equation}
   (Here, the bar denotes complex conjugation: $\overline{z}=(\overline{z}_{1},\ldots,\overline{z}_{d})$).
   It follows from (\ref{eq3.1}) that as the coefficients of the denominator $q(z)$, we can always choose real numbers ($\overline{q(\overline{z})}=q(z)$).
   Then the coefficients of the numerator $P(z)$ are Hermitian $m\times m$ matrices.

   It turn out that every rational Pick function in $d$ variables can be obtained from the Cayley inner function in $d+1$ variables:

  \begin{theorem}\label{th3.1}
      Let $f(z)=P(z)/q(z)$ be a rational $\mathbb{C}^{m\times m}$-valued Pick function in $d$ variables,
      where $P(z)$, $q(z)$ are coprime polynomials. If
          \begin{multline}\label{eq3.2}
          \quad P_{1}(z)=[\,P(z)-P(\overline{z})^{\ast}\,]/2i,\quad P_{2}(z)=[\,P(z)+P(\overline{z})^{\ast}\,]/2,\\
          q_{1}(z)=[\,q(z)-\overline{q(\overline{z})}\,]/2i,\quad q_{2}(z)=[\,q(z)+\overline{q(\overline{z})}\,]/2,\qquad\enskip
          \end{multline}
      then
             \begin{equation}\label{eq3.3}
             g(z_{0},z)=\frac{z_{0}P_{1}(z)+P_{2}(z)}{z_{0}q_{1}(z)+q_{2}(z)}
             \end{equation}
      is the rational Cayley inner function in $(d+1)$ variables of the Pick class. Moreover,
               \begin{equation}\label{eq3.4}
               g(i,z)=f(z).
               \end{equation}
   \end{theorem}
   For proof this, we need several known statements.

    We say that a multivariate polynomial $p(z)$ with complex coefficients is
    \emph{stable}\footnote{This terminology may differ from the designations of other authors.}
    if it is nonzero whenever $z\in\Pi^{d}$. A stable polynomial with real coefficients will be called \emph{a real stable}.
    The ring of the polynomials with real coefficients is denoted by $\mathbb{R}[z_{1},\ldots,z_{d}]$.
    \medskip

   \textbf{Lemma 3.1} (\cite{uj04}, Corollary 5.5, \cite{uj05}, Lemma 2.2).
        \emph{Let $p(z)+iq(z)\neq0$, where $p,q\in\mathbb{R}[z_{1},\ldots,z_{d}]$ and let $z_{d+1}$ be a new indeterminate.
         Then the following are equivalent.}

         (a) $p(z)+iq(z)$ \emph{is stable,}

         (b) $p(z)+z_{d+1}q(z)$ \emph{is real stable,}

         (c) \emph{all nonzero polynomials in the pencil
                $$
                \{\alpha p(z)+\beta q(z)\mid \alpha,\beta\in\mathbb{R}\}
                $$
              are real stable,}

         (d) $\mbox{Im}\,p(z)/q(z)\geq0$ \emph{whenever $\emph{Im}\,z_{i}>0$ for all $1\leq i\leq d$.}
     \medskip

  \emph{Proof of Theorem} \ref{th3.1}.
      Since $f(z)=P(z)/q(z)\in\mathcal{P}_{d}^{m\times m}$, we see that
                       $$
                        f_{\eta}(z)=\frac{\eta P(z)\eta^{\ast}}{q(z)}
                       $$
      is the scalar function of the Pick class for every row vector $\eta\in\mathbb{C}^{m}$.
      $f_{\eta}(z)+\alpha$, where $\alpha\in\mathbb{R}$, is also function of the Pick class.
      Then its numerator $\alpha q(z)+\eta P(z)\eta^{\ast}$ is a stable polynomial with complex coefficients.
      Using (\ref{eq3.2}), we represent it the form
         $$
         \alpha q(z)+\eta P(z)\eta^{\ast}=i(\alpha q_{1}(z)+\eta P_{1}(z)\eta^{\ast})+(\alpha q_{2}(z)+\eta P_{2}(z)\eta^{\ast}),
         $$
      where $(\alpha q_{1}(z)+\eta P_{1}(z)\eta^{\ast}),\,(\alpha q_{2}(z)+\eta P_{2}(z)\eta^{\ast})\in\mathbb{R}[z_{1},\ldots,z_{d}]$.

      Let $z_{0}$ be a new indeterminate. By Lemma 3.1, for every $\eta\in\mathbb{C}^{m}$
           $$
           \alpha (z_{0}q_{1}(z)+q_{2}(z))+\eta(z_{0}P_{1}(z)+P_{2}(z))\eta^{\ast}
           $$
      is the real stable polynomial in $(d+1)$ variables. Since $\eta\in\mathbb{C}^{m}$ is an arbitrary vector,
      we see that for every $\alpha,\beta\in\mathbb{R}$ all nonzero polynomials in the pencil
              $$
              \alpha(z_{0}q_{1}(z)+q_{2}(z))+\beta(\eta(z_{0}P_{1}(z)+P_{2}(z))\eta^{\ast})
              $$
      are real stable. Then
                  $$
                  (z_{0}q_{1}(z)+q_{2}(z))+z_{d+1}(\eta(z_{0}P_{1}(z)+P_{2}(z))\eta^{\ast})
                  $$
      is the real stable polynomial in $(d+2)$ variables. Therefore, for every $\eta\in\mathbb{C}^{m}$
                     $$
                     \mbox{Im}\,\frac{\eta(z_{0}P_{1}(z)+P_{2}(z))\eta^{\ast}}{z_{0}q_{1}(z)+q_{2}(z)}\geq0,\,
                     \mbox{whenever}\,\, \mbox{Im}\,z_{i}>0\,\, \mbox{for all}\,\, 0\leq i\leq d.
                     $$
      From this
           $$
           g(z_{0},z)=\frac{z_{0}P_{1}(z)+P_{2}(z)}{z_{0}q_{1}(z)+q_{2}(z)}\in\mathcal{P}_{d+1}^{m\times m}.
           $$
      For real values of variables, the values of $g(x_{0},x)$ are Hermitian matrices.
      Then $g(z_{0},z)\in\mathcal{IP}_{d+1}^{m\times m}$. Identity (\ref{eq3.4}) is obvious.
      \qed
      \medskip

    We say that a matrix-valued polynomial $F(z)$ is \emph{positive semidefinite} or \emph{PSD}
    if $F(x)\geq0$ for all $x\in\mathbb{R}^{d}$.
    \medskip

    For Cayley inner Pick functions, the following statement holds:
        \begin{proposition}\label{pr3.1}
            If $f(z) = P(z)/q(z)\in \mathcal{IP}_{d}^{m\times m}$, then partial Wronskians
                \begin{equation}\label{eq3.5}
                W_{k}[q,\,P]=q(z)\frac{\partial P}{\partial z_{k}}(z)-P(z)\frac{\partial q}{\partial z_{k}}(z),\quad k=1,\ldots,d
                \end{equation}
            are PSD polynomials.
      \end{proposition}

   \begin{proof}
      Suppose $k=1$. Let us $\varphi(\zeta)=f(\zeta,\widehat{x})$. If $\widehat{x}=(x_{2},\ldots,x_{d})\in\mathbb{R}^{d-1}$, then
          $\mbox{Im}\,\varphi (\zeta)\geq 0\enskip \mbox{for}\enskip \mbox{Im}\,\zeta>0$
      and $\mbox{Im}\,\varphi (\zeta)=0\enskip \mbox{for}\enskip\mbox{Im}\,\zeta=0$.
      Hence inequality $d\varphi(\zeta)/d\zeta \mid_{\zeta\in\mathbb{R}}\geq 0$ holds. From this
      \medskip

      \qquad \qquad $W_{z_{1}}[\,q,P\,](x)=q(x)^{2}d\varphi(\zeta)/ d\zeta|_{\zeta=x_{1}}\geq0, \quad x\in\mathbb{R}^{d}$.
   \end{proof}
       For multi-affine functions, the nonnegativity of the Wronskians is also a sufficient condition.
      \begin{theorem}\label{th3.2}
        A rational multi-affine matrix-valued function
             $f(z)=P(z)/q(z)$ belongs to the class $\mathcal{IP}_{d}^{m\times m}$ if and only if all Wronskians
                 $$
                 W_{k}[\,q,P\,]=q(z)\frac{\partial P}{\partial z_{k}}(z)-P(z)\frac{\partial q}{\partial z_{k}}(z),\quad k=1,\ldots,d
                 $$
       are matrix-valued PSD polynomials.
      \end{theorem}
         To prove this theorem, we need the following statement.
   \begin{theorem}\label{th3.3}
        Let $f(z)$ be a rational matrix-valued function.
        Let us assume that for every  $j=1,2,\ldots,d$ and for every real $x_{1},x_{2},\ldots,x_{d}$, the functions
            $$
            f_{j}(z)=f(x_{1},\ldots,x_{j-1},z_{j},x_{j+1,}\ldots,x_{d})
            $$
        satisfy the conditions
               $$
               (f_{j}-f_{j}^{\ast})/2i\geq0\quad \mbox{for}\quad \emph{Im}\,z_{j}>0.
               $$
        Then the inequality
                  $$
                  (f(z)-f(z)^{\ast})/2i\geq0\quad\mbox{for}\quad z\in \Pi^{d}
                  $$
        holds.
   \end{theorem}
          \begin{proof}
             For row vector $\xi$ does not depended on $z$, let us consider the function $f_{\xi}(z)=\xi f(z)\xi^{\ast}$.
             It is clear that
                  $$
                  \mbox{Im}\,f_{\xi}(z)=\xi[\,(f(z)-f(z)^{\ast})/2 i\,]\xi^{\ast}.
                  $$
             Therefore, if the inequality
                    $$
                    \mbox{Im}\,f_{\xi}(z)\geq0\quad\mbox{for}\quad z\in \Pi^{d}
                    $$
             holds for every $\xi$, then the inequality
                       $$
                       (f(z)-f(z)^{\ast})/2i\geq0\quad\mbox{for}\quad z\in \Pi^{d}
                       $$
             holds us well. Therefore, it is enough to consider scalar functions only.
             Let a scalar function $f(z)$ satisfy the assumptions of the theorem.
             We consider rational function
                          \begin{equation}\label{eq3.6}
                           u(\zeta_{1},\ldots,\zeta_{d})=\frac{f\left(i\frac{1+\zeta_{1}}{1-\zeta_{1}},\ldots,
                           i\frac{1+\zeta_{d}}{1-\zeta_{d}}\right)-i}{f\left(i\frac{1+\zeta_{1}}{1-\zeta_{1}},
                          \ldots,i\frac{1+\zeta_{d}}{1-\zeta_{d}}\right)+i}.
                          \end{equation}
             Let $\mathbb{T}$ be the unit circle. it is clear, that for every $j=1,2,\ldots,d$ and for every
             $t_{1}\in\mathbb{T},\ldots,t_{j-1}\in\mathbb{T},t_{j+1}\in\mathbb{T},\ldots,t_{d}\in\mathbb{T}$, the function
             $u(t_{1},\ldots,t_{j-1},\zeta_{j},t_{j+1},\ldots,t_{d})$ is golomorphic for $|\zeta_{j}|<1$ and satisfies the inequality
                   \begin{equation}\label{eq3.7}
                   |u(t_{1},\ldots,t_{j-1},\zeta_{j},t_{j+1},\ldots,t_{d})|<1\quad\mbox{for}\quad |\zeta_{j}|<1.
                   \end{equation}
             Let us show that the function $u(\zeta_{1},\ldots,\zeta_{d})$ is golomorphic in polydisk
             $\mathbb{D}^{d}=\{\zeta\in\mathbb{C}^{d}\mid |\zeta_{1}|<1,\ldots,|\zeta_{d}|<1\}$ and satisfies the inequality
                       $$
                       |u(\zeta_{1},\ldots,\zeta_{d})|<1\quad\mbox{for}\quad\zeta\in\mathbb{D}^{d}.
                       $$
             To prove this, we consider the Fourier coefficients $\widehat{u}(k_{1},\ldots,k_{d})$ of the function
             $u(t_{1},\ldots,t_{d})$ considered on the torus $\mathbb{T}^{d}=\{t\in\mathbb{C}^{d}\mid |t_{1}|=1,\ldots,|t_{d}|=1\}$:
                 $$
                 \widehat{u}(k_{1},\ldots,k_{d})=\int_{\mathbb{T}^{d}}u(t_{1},\ldots,t_{d})t_{1}^{-k_{1}}\cdots t_{d}^{-k_{d}}m(dt_{1})\cdots m(dt_{d}).
                 $$
             ($m(dt)$ is one-dimensional normalized Lebesgue measure). The function $u$ is contractive on $\mathbb{T}^{d}$: $|u(t)|\leq1$ for $t\in\mathbb{T}^{d}$.
             Therefore, its Fourier coefficients  $\widehat{u}(k_{1},\ldots,k_{d})$ exist. If $k_{j}<0$ at least for one $j=1,2,\ldots,d$,
             then $\widehat{u}(k_{1},\ldots,k_{d})=0$. Indeed, for definiteness, let $k_{1}<0$. Then
                     \begin{multline}\label{eq3.8}
                     \widehat{u}(k_{1},\ldots,k_{d})=\int_{\mathbb{T}^{d-1}}t_{2}^{-k_{2}}\cdots t_{d}^{-k_{d}}m(dt_{2})\cdots m(dt_{d})\times \\
                     \int_{\mathbb{T}}u(t_{1},\ldots,t_{d})t_{1}^{-k_{1}}m(dt_{1}).
                     \end{multline}
             By condition (\ref{eq3.7}), the inner integral in (\ref{eq3.8}) vanishes. Therefore, the Fourier coefficients $\widehat{u}(k_{1},\ldots,k_{d})$
             determine the function
                 $$
                 g(\zeta_{1},\ldots,\zeta_{d})=\sum_{\forall k}\widehat{u}(k_{1},\ldots,k_{d})\zeta_{1}^{k_{1}}\cdots\zeta_{d}^{k_{d}},
                 $$
             which is golomorphic in $\mathbb{D}^{d}$. On the other hand, denoting $\zeta_{j}=r_{j}t_{j}\enskip(r_{j}\geq0,\enskip t_{j}\in\mathbb{T})$,
             we obtain
                   $$
                   g(\zeta_{1},\ldots,\zeta_{d})=\sum_{\forall k}r_{1}^{|k_{1}|}\cdots r_{d}^{|k_{d}|}\widehat{u}(k_{1},\ldots,k_{d})t_{1}^{k_{1}}\cdots t_{d}^{k_{d}}.
                   $$
             Therefore, $g(\zeta_{1},\ldots,\zeta_{d})$ is the convolution of the function $u(t_{1},\ldots,t_{d})$ and the Poisson kernel
                      $$
                      P(r,t)=\sum_{\forall k}r_{1}^{|k_{1}|}\cdots r_{d}^{|k_{d}|}t_{1}^{k_{1}}\cdots t_{d}^{k_{d}}.
                      $$
             Since $|u(t_{1},\ldots,t_{d}|\leq1$ on $\mathbb{T}^{d}$ and
                         $$
                         \int_{\mathbb{T}^{d}}P(r_{1},\ldots,r_{d},t_{1},\ldots,t_{d})m(dt_{1})\cdots m(dt_{d})=1,
                         $$
             we have
                  $$
                  |g(\zeta_{1},\ldots,\zeta_{d}|\leq1\quad\mbox{for}\quad \zeta\in\mathbb{D}^{d}.
                  $$
             Since the Poisson kernel is approximate identity,
                     $$
                     \lim_{r\rightarrow1-0}g(r_{1}t_{1},\ldots,r_{d}t_{d})=u(t_{1},\ldots,t_{d})
                     $$
             in every point $(t_{1},\ldots,t_{d})\in\mathbb{T}^{d}$ where the function $u$ is continuous.
             By the uniqueness theorem,
                        $$
                        u(\zeta)=g(\zeta)\quad\mbox{for}\quad \zeta\in\mathbb{D}^{d}.
                        $$
             Returning to $f(z_{1},\ldots,z_{d})$ by means of the transformation that is inverse
             to the transformation (\ref{eq3.6}) we obtain the statement of theorem.
             \end{proof}

     \noindent
     \emph{Proof of Theorem} \ref{th3.2}.
          The necessity is proved in Proposition \ref{pr3.1}. Let us prove the sufficiency. Since $f(z)$ is multi-affine, we see that
              $$
              f(z)=\frac{z_{k}P_{1}(\widehat{z})+P_{2}(\widehat{z})}{z_{k}q_{1}(\widehat{z})+q_{2}(\widehat{z})},
              \quad \widehat{z}=(z_{1},\ldots,z_{k-1},z_{k+1},\ldots,z_{d}).
              $$
          From here
                 \begin{equation}\label{eq3.9}
                 \mbox{Im}\,f(z_{k},\widehat{x})=\mbox{Im}\,z_{k}\,\frac{W_{k}[\,q,P\,](\widehat{x})}{| z_{k}q_{1}(\widehat{x})+q_{2}(\widehat{x})| ^{2}},
                  \quad \widehat{x}\in \mathbb{R}^{d-1}.
                 \end{equation}
          Hence $\mbox{Im}\,f(z_{k},\widehat{x})\geq 0, \mbox{Im} \,z_{k}>0$ for each $k=1,\ldots,d$ (for any other real variables).
          By Theorem \ref{th3.3}, $\mbox{Im}\,f(z)\geq 0$ for $z\in \Pi^{d}$. If $\mbox{Im}\,z_{k}=0$, $k=1,\ldots,d$, then
          $\mbox{Im}\,f(z_{k},\widehat{z})=0$, that is, $f(z)\in\mathcal{IP}_{d}^{m\times m}$.
     \qed

   \section{PSD Polynomials and Sums of Squares}\label{s:4}

   Recall that a matrix-valued polynomial $F(z)$ is called \emph{a positive semidefinite} or \emph{PSD}
   if its values $F(x)$ are positive semidefinite Hermitian matrices: $F(x)\geq0$ for every $x\in\mathbb{R}^{d}$.

   We say that a $\mathbb{C}^{m\times m}$-valued PSD polynomial $F(z)$ is \emph{a sum of squares} or \emph{SOS} if
           \begin {equation}\label{eq4.1}
           F(z)=H(z)H(\overline{z})^{\ast},
           \end{equation}
   where $H(z)$ is some $\mathbb{C}^{m\times k}$-valued polynomial.

     Not every PSD polynomial is an SOS. An amazing fact is that the PSD polynomials (\ref{eq3.5})
     associated with rational Cayley inner functions of the Pick class are SOS polynomials:
        \begin{theorem}\label{th4.1}\emph{\textbf{(Sum-of-Squares Theorem)}}.
           If $f(z)=P(z)/q(z)\in\mathcal{IP}_{d}^{m\times m}$, then
               \begin{equation}\label{eq4.2}
                W_{k}[q,P]=q(z)\frac{\partial P}{\partial z_{k}}(z)-P(z)\frac{\partial q}{\partial z_{k}}(z),\quad k=1,\ldots,d
               \end{equation}
           are matrix-valued SOS polynomials.
        \end{theorem}

     For the convenience of the reader, the proof of this theorem is considered in Appendix.

   \section{Generalization of Darlington's Theorem}\label{s:5}

     The application of Darlington's theorem to obtain representations of matrix functions in one variable of some classes was considered in \cite{uj20}.
     In \cite{uj18}, an analogue of Darlington's method was used for a positive real functions of several variables.
     This method allows us to reduce the question of representing a function in $d$ variables
     to the question of representing a similar function in $d-1$ variables.

     For multi-affine Cayley inner functions of the Pick class, the following generalization of Darlington theorem holds.

        \begin{theorem}\label{th5.1} \enskip
          Let a multi-affine function $f(z_{0},z)\in\mathcal{IP}_{d+1}^{m\times m}$ be represented as
              $$
              f(z_{0},z)=\frac{P(z_{0},z)}{q(z_{0},z)}=\frac{z_{0}P_{1}(z)+P_{2}(z)}{z_{0}q_{1}(z)+q_{2}(z)},\quad z\in\mathbb{C}^{d},
              $$
          and let
                 $$
                  W_{z_{0}}[\,q,P\,]=P_{1}(z)q_{2}(z)-P_{2}(z)q_{1}(z)=\Phi_{0}(z){\Phi}_{0}(\overline z)^{\ast}
                 $$
          be an SOS polynomial. Suppose $\Phi_{0}(z)$ has the size $m\times r$. Then

          \textbf{\emph{(i)}} If $q_{1}(z)\neq0$, then
               \begin{equation}\label{eq5.1}
                  g(z)=
                    \begin{pmatrix}
                     g_{11}(z) & g_{12}(z) \\
                     g_{21}(z) & g_{22}(z)
                    \end{pmatrix}=
                         \frac{1}{q_{1}(z)}
                                \begin{pmatrix}
                                          P_{1}(z)            &  \Phi_{0}(z) \\
                                \Phi_{0}(\overline{z})^{\ast} & q_{2}(z)I_{r}
                                \end{pmatrix}
               \end{equation}
          is a multi-affine function of class $\mathcal{IP}_{d}^{(m+r)\times (m+r)}$ and
                    \begin{equation}\label{eq5.2}
                    f(z_{0},z)=g_{11}(z)-g_{12}(z)(g_{22}(z)+z_{0}I_{r})^{-1}g_{21}(z).
                    \end{equation}

          \textbf{\emph{(ii)}} If $q_{2}(z)\neq0$, then
               \begin{equation}\label{eq5.3}
                   g(z)=
                         \begin{pmatrix}
                         g_{11}(z) & g_{12}(z) \\
                         g_{21}(z) & g_{22}(z)
                         \end{pmatrix}=
                              \frac{1}{q_{2}(z)}
                                     \begin{pmatrix}
                                               P_{2}(z)            &  \Phi_{0}(z) \\
                                     \Phi_{0}(\overline{z})^{\ast} & -q_{1}(z)I_{r}
                                     \end{pmatrix}
               \end{equation}
          is a multi-affine function of class $\mathcal{IP}_{d}^{(m+r)\times (m+r)}$ and
                    \begin{equation}\label{eq5.4}
                    f(z_{0},z)=g_{11}(z)-g_{12}(z)(g_{22}(z)-z_{0}^{-1}I_{r})^{-1}g_{21}(z).
                    \end{equation}
       \end{theorem}

        \noindent
        \emph{Proof}.
        Both cases are proved in similar way. Lt us prove, for example, \textbf{(ii)}.
        Representation (\ref{eq5.4}) follows from the obvious identity
              \begin{equation}\label{eq5.5}
               f(z_{0},z)=\frac{z_{0}P_{1}(z)+P_{2}(z)}{z_{0}q_{1}(z)+q_{2}(z)}=\frac{P_{2}}{q_{2}}-
              \frac{\Phi_{0}(z)\Phi_{0}(\overline{z})^{\ast}}{q_{2}^{2}(-q_{1}/q_{2}-z_{0}^{-1})}.
              \end{equation}
        The multi-affinity of $g(z)$ is obvious. Let us prove $g(z)\in\mathcal{IP}_{d}^{(m+r)\times (m+r)}$.
        By Theorem \ref{th3.2}, it suffices to prove $F_{k}(z)=q_{2}^{2} \partial g(z)/\partial z_{k},\, k=1,\ldots,d$
        are PSD polynomials. $f(z_{0},z)$ is multi-affine. Then
                   \begin{equation}\label{eq5.5}
                   f(z_{0},z)=\frac{P(z_{0},z)}{q(z_{0},z)}=
                   \frac{z_{0}z_{k}\widehat{P}_{1}+z_{0}\widehat{P}_{2}+z_{k}\widehat{P}_{3}+\widehat{P}_{4}}
                   {z_{0}z_{k}\widehat{q}_{1}+z_{0}\widehat{q}_{2}+z_{k}\widehat{q}_{3}+\widehat{q}_{4}}.
                   \end{equation}
        From (\ref{eq5.5}) and (\ref{eq5.3}) we get
             \begin{multline}\label{eq5.6}
              W_{z_{0}}[q,P]=\Phi_{0}(z)\Phi_{0}(\overline{z})^{\ast}=
              z_{k}^{2}(\widehat{q}_{3}\widehat{P}_{1}-\widehat{q}_{1}\widehat{P}_{3})+\\
             z_{k}(\widehat{q}_{4}\widehat{P}_{1}-\widehat{q}_{1}\widehat{P}_{4}
             +\widehat{q}_{3}\widehat{P}_{2}-\widehat{q}_{2}\widehat{P}_{3})
             +(\widehat{q}_{4}\widehat{P}_{2}-\widehat{q}_{2}\widehat{P}_{4}),
             \end{multline}
                 \begin{equation}\label{eq5.7}
                 F_{k}=q_{2}^{2}\frac{\partial g(z)}{\partial z_{k}}=
                     \begin{pmatrix}
                     \widehat{P}_{3}\widehat{q}_{4}-\widehat{P}_{4}\widehat{q}_{3} &             \Phi_{k}(z) \\
                     \Phi_{k}(\overline{z})^{\ast}                      & (\widehat{q}_{2}\widehat{q}_{3}-\widehat{q}_{1}\widehat{q}_{4})I_{r}
                     \end{pmatrix},
                 \end{equation}
        where $\Phi_{k}(z)=(z_{k}\widehat{q}_{3}+\widehat{q}_{4})\partial\Phi_{0}/\partial z_{k}-\widehat{q}_{3}\Phi_{0},\,k=1,\ldots,d$.
        Note that $\Phi_{k}(z)$ is actually independent of $z_{0}$ and $z_{k}$.

        The diagonal elements $f_{ii}(z_{0},z)$ of a rational matrix-valued function $f(z_{0},z)\in\mathcal{IP}_{d}^{m\times m}$
        satisfy the condition $\overline{f_{ii}(\overline{z}_{0},\overline{z})}=f_{ii}(z_{0},z)$.
        Since $f_{ii}=p_{ii}/q$, we see that the coefficients of the polynomial $q(z_{0},z)$ can be considered real numbers.
        Note the identity
             \begin{equation}\label{eq5.8}
             (\widehat{P}_{3}\widehat{q}_{4}-\widehat{P}_{4}\widehat{q}_{3})=\Phi_{k}(z)(\widehat{q}_{2}\widehat{q}_{3}-\widehat{q}_{1}
             \widehat{q}_{4})^{-1}\Phi_{k}(\overline{z})^{\ast},
             \end{equation}
        or, equivalently,
                \begin{equation}\label{eq5.81}
                (\widehat{P}_{3}\widehat{q}_{4}-\widehat{P}_{4}\widehat{q}_{3})(\widehat{q}_{2}\widehat{q}_{3}-\widehat{q}_{1}
                \widehat{q}_{4})=\Phi_{k}(z)\Phi_{k}(\overline{z})^{\ast}.
                \end{equation}
        Indeed, since the coefficients of the polynomials $\widehat{q_{i}}(z)$ are real, then
            \begin{multline}\label{eq5.9}
            \Phi_{k}(z)\Phi_{k}(\overline{z})^{\ast}=(z_{k}\widehat{q}_{3}+\widehat{q}_{4})^{2}\frac{\partial \Phi_{0}(z)}{\partial z_{k}}
            \frac{\partial \Phi_{0}(\overline{z})^{\ast}}{\partial z_{k}}-\\
            (z_{k}\widehat{q}_{3}^{2}+\widehat{q}_{3}\widehat{q}_{4})\left(\frac{\partial\Phi_{0}(z)}{\partial z_{k}}\Phi_{0}(\overline{z})^{\ast}
            +\Phi_{0}(z)\frac{\partial\Phi_{0}(\overline{z})^{\ast}}{\partial z_{k}}\right)
            +\widehat{q}_{3}^{2}\Phi_{0}(z)\Phi_{0}(\overline{z})^{\ast}.
            \end{multline}
        $\Phi_{0}(\widehat{z})$ is a multi-affine polynomial. Differentiating (\ref{eq5.6}) in $z_{k}$ and substituting
        the obtained expressions into (\ref{eq5.9}), we get (\ref{eq5.81}).

        Since
            $$
            h=\left.\frac{q}{\partial q/\partial z_{0}}\right|_{z_{0}=0}=
            \frac{z_{k}\widehat{q}_{3}+\widehat{q}_{4}}{z_{k}\widehat{q}_{1}+\widehat{q}_{2}}\in\mathcal{IP}_{d},
            $$
        we see that $Q(z)=\widehat{q}_{2}(z)\widehat{q}_{3}(z)-\widehat{q}_{1}(z)\widehat{q}_{4}(z)$ is PSD polynomial.
        For $x\in \mathbb{R}^{d-1}$, from (\ref{eq5.7}), (\ref{eq5.8}) we obtain
                 $$
                 F_{k}(x)=
                     \begin{pmatrix}
                     I_{m}  & \Phi_{k}(x)Q^{-1}(x)\\
                        0   & I_{r}
                     \end{pmatrix}
                          \begin{pmatrix}
                            0 &     0 \\
                            0 & Q(x)I_{r}
                          \end{pmatrix}
                               \begin{pmatrix}I_{m}      & 0 \\
                                Q^{-1}(x)\Phi_{k}(x)^{T} & I_{r}
                               \end{pmatrix}\geq 0,
                 $$
        i.e., $F_{k}(z)$, $k=1,\ldots,d$ are the PSD polynomials.
        By Theorem \ref{th3.2}, the function $g(z)$ (\ref{eq5.3}) belongs to the class $\mathcal{IP}_{d}^{(m+r)\times(m+r)}$
        \qed

    \section{Superposition of Coefficient Matrix\\ of Fractional Linear Transformations}\label{s:6}
      We will need a somewhat unusual superposition of fractional linear transformations of the form
          \begin{equation}\label{eq6.1}
          W=g_{11}-g_{12}(g_{22}+Z)^{-1}g_{21}.
          \end{equation}
      The unusual thing is that the argument of transformation (\ref{eq6.1}) is assumed to be the matrix
              $$
               \begin{pmatrix}
               g_{11} & g_{12} \\
               g_{21} & g_{22}
               \end{pmatrix},
              $$
      and the matrix $Z$ plays the role of parameter defining transformation (\ref{eq6.1}).
      When considering the superposition of such transformations, it is convenient to use the following
      notation for transformation (\ref{eq6.1}):
                    \begin{equation}\label{eq6.2}
                    W=g_{11}-g_{12}(g_{22}+Z)^{-1}g_{21}:=
                         \left(
                           \begin{array}{c|c}
                           g_{11} & g_{12} \\
                           \hline
                           g_{21} & g_{22}+Z
                           \end{array}
                         \right).
                    \end{equation}

      The following statement holds.
      \begin{proposition}\label{pr6.1}
            If
              \begin{equation}\label{eq6.3}
                 W=
                   \left(
                   \begin{array}{c|c}
                   g_{11} & g_{12} \\
                   \hline
                   g_{21} & g_{22}+Z_{1}
                   \end{array}
                   \right)
              \end{equation}
           and
                  \begin{equation}\label{eq6.4}
                    \begin{pmatrix}
                    g_{11} & g_{12}\\
                    g_{21} & g_{22}
                    \end{pmatrix}=
                       \left(
                       \begin{array}{cc|c}
                       a_{11} & a_{12} & a_{13} \\
                       a_{21} & a_{22} & a_{23} \\
                       \hline
                       a_{31} & a_{32} & a_{33}+Z_{2}
                       \end{array}
                       \right),
                  \end{equation}
          then
             \begin{equation}\label{eq6.5}
               W=
             \left(
             \begin{array}{c|cc}
             a_{11} &    a_{12}    &   a_{13}   \\
             \hline
             a_{21} & a_{22}+Z_{1} &   a_{23}    \\
             a_{31} &    a_{32}    & a_{33}+Z_{2}
             \end{array}
             \right).
          \end{equation}
     \end{proposition}

     \begin{proof}
     It follows from (\ref{eq6.4}) and (\ref{eq6.3}) that the matrices $S=(a_{33}+Z_{2})$ and $g_{22}+Z_{1}=(a_{22}+Z_{1})-a_{23}(a_{33}+Z_{2})^{-1}a_{32}=$
     are invertible. By the Frobenius formula \cite{uj08}, there exists
             \begin{multline}\label{eq6.6}
                \begin{pmatrix}
                a_{22}+Z_{1} &    a_{23}\\
                    a_{32}   & a_{33}+Z_{2}
                \end{pmatrix}^{-1}=\\
                       \begin{pmatrix}
                               (g_{22}+Z_{1})^{-1}       &          -(g_{22}+Z_{1})^{-1}a_{23}S^{-1}          \\
                       -S^{-1}a_{32}(g_{22}+Z_{1})^{-1}  & S^{-1}+S^{-1}a_{32}(g_{22}+Z_{1})^{-1}a_{23}S^{-1}
                       \end{pmatrix}.
             \end{multline}
     For the right-hand side of (\ref{eq6.5}) we obtain
                \begin{multline}\label{eq6.7}
                a_{11}-(a_{12},\,a_{13})\times\\
                     \begin{pmatrix}
                        (g_{22}+Z_{1})^{-1}            &          -(g_{22}+Z_{1})^{-1}a_{23}S^{-1}  \\
                     -S^{-1}a_{32}(g_{22}+Z_{1})^{-1}  & S^{-1}+S^{-1}a_{32}(g_{22}+Z_{1})^{-1}a_{23}S^{-1}
                     \end{pmatrix}
                           \begin{pmatrix}
                           a_{21} \\
                           a_{31}
                           \end{pmatrix}=\\
                a_{11}-a_{12}(g_{22}+Z_{1})^{-1}a_{21}+a_{13}S^{-1}a_{32}(g_{22}+Z_{1})^{-1}a_{21}+\\
                a_{12}(g_{22}+Z_{1})^{-1}a_{23}S^{-1}a_{31}-\\
                     a_{13}S^{-1}a_{31}-a_{13}S^{-1}a_{32}(g_{22}+Z_{1})^{-1}a_{23}S^{-1}a_{31}=\\
                        (a_{11}-a_{13}S^{-1}a_{31})-(a_{12}-a_{13}S^{-1}a_{32})(g_{22}+Z_{1})^{-1}a_{21}+\\
                            (a_{11}-a_{13}S^{-1}a_{31})(g_{22}+Z_{1})^{-1}a_{23}S^{-1}a_{31}=\\
                (a_{11}-a_{13}S^{-1}a_{31})-(a_{12}-a_{13}S^{-1}a_{32})(g_{22}+Z_{1})^{-1}(a_{21}-a_{23}S^{-1}a_{31})=\\
                     g_{11}=g_{12}(g_{22}+Z_{1})^{-1}g_{21}=W,\qquad\quad
                \end{multline}
     which was required.
     \end{proof}

   \section{Representation of Cayley inner functions\\ of the Pick class}\label{s:7}

   For Cayley inner rational functions of the Pick class, the following statement holds:

   \begin{theorem}\label{th7.1}
       Let $f$ be a $\mathbb{C}^{m\times m}$--valued function of $d$ complex variables. The following statements are equivalent.
       \medskip

       \textbf{\emph{(0)}} $f(z)$ is a rational Cayley inner function of the Pick class.
     \medskip

     \textbf{\emph{(1)}} There exist $n, n_{1},\ldots,n_{d}\in\mathbb{Z}_{+}$, $ n_{0}\geq0$, $n=n_{0}+n_{1}+\cdots+n_{d}$ and the $z$-independent
     Hermitian $(m+n)\times (m+n)$ matrix
                $$
                B=
                  \begin{pmatrix}
                    B_{11} & B_{12} \\
                    B_{21} & B_{22}
                  \end{pmatrix},\quad B=B^{\ast}
                $$
     with block $B_{22}$ of size $n\times n$ such that
              \begin{equation}\label{eq7.1}
                  f(z)=B_{11}-B_{12}(B_{22}+Z_{n})^{-1}B_{21}
              \end{equation}
     where $Z_{n}=\emph{diag}\{0_{n_{0}},z_{1}I_{n_{1}},\ldots,z_{d}I_{n_{d}}\}$.
     \end{theorem}

     \begin{proof}
        \textbf{(0)} $\Rightarrow$ \textbf{(1)}.
        Let $P(z),\,q(z)$ be coprime polynomials and $f(z)=P(z)/q(z)\in\mathcal{IP}_{d}^{m\times m}$. If
            $$
            \max \{\mbox{deg}_{z_{j}}P(z),\mbox{deg}_{z_{j}}q(z)\}=k_{j},\quad j=1,\ldots,d,
            $$
        then, by Theorem \ref{th2.1}, a multi-affine function
               \begin{equation}\label{eq7.2}
               h_{0}(\zeta_{1},\ldots,\zeta_{k})=\mathbf{D}_{z_{1}}^{k_{1}}\cdots\mathbf{D}_{z_{d}}^{k_{d}}[f(z_{1},\ldots,z_{d})],
               \end{equation}
        where $k=k_{1}+\cdots+k_{d}$, belongs to the class $\mathcal{IP}_{k}^{m\times m}$.
        Note that the set of variables $\{\zeta_{1},\ldots,\zeta_{k}\}$ of the function (\ref{eq7.2}) is the union of $d$ subsets.
        By identifying all variables of every $i$-th subset with the variable $z_{i}$, we obtain
        the original function $f(z_{1},\ldots,z_{d})$.
        \medskip

        We will construct the representation of the function (\ref{eq7.2}) step by step.
        At the $j$-th step, we will obtain a multi-affine Cayley inner matrix-valued function $h_{j}(\zeta)$
        of the Pick class, depending on the $(k-j)$ variables.
        \medskip

        \noindent
        \textbf{\emph{Step} 1}.
             Since $h_{0}(\zeta)$ is multi-affine, we see that
                  \begin{equation}\label{eq7.3}
                  h_{0}(\zeta_{1},\ldots,\zeta_{k})=\frac{\zeta_{1}P_{1}(\zeta)+P_{2}(\zeta)}{\zeta_{1}q_{1}(\zeta)+q_{2}(\zeta)},
                  \end{equation}
             where $P_{1},P_{2},q_{1},q_{2}$ are independent of $\zeta_{1}$. By Proposition \ref{pr3.1}, the Wronskian
                      \begin{equation}\label{eq7.4}
                      W_{\zeta_{1}}[q,P]=P_{1}(\zeta)q_{2}(\zeta)-P_{2}(\zeta)q_{1}(\zeta)
                      \end{equation}
             is the PSD polynomial.

             By Sum-of-Squares Theorem, there exists $\mathbb{C}^{m\times r_{1}}$-valued polynomial $\Phi_{1}(\zeta)$ such that
                      $$
                      W_{\zeta_{1}}[q,P]=\Phi_{1}(\zeta)\Phi_{1}(\overline{\zeta})^{\ast}.
                      $$
             The polynomial $\Phi_{1}(\zeta)$ does not depend on $\zeta_{1}$ and has degree at most 1 for every variable.

             For the denominator $\zeta_{1}q_{1}(\zeta)+q_{2}(\zeta)$ in (\ref{eq7.3}) there are 2 possibilities:

             \textbf{(i)} $q_{1}(\zeta)\neq0$.
                By Theorem \ref{th5.1},
                      \begin{equation}\label{eq7.5}
                      h_{0}(\zeta)=g_{11}-g_{12}(g_{22}+\zeta_{1}I_{r_{1}})^{-1}g_{21}=
                          \left(
                          \begin{array}{c|c}
                          g_{11}   &   g_{12}\\
                          \hline
                          g_{21}   & g_{22}+\zeta_{1}I_{r_{1}}
                          \end{array}
                         \right),
                      \end{equation}
                where multi-affine function
                          \begin{equation}\label{eq7.6}
                          g_{1}(\zeta)=
                               \begin{pmatrix}
                               g_{11} & g_{12} \\
                               g_{21} & g_{22}
                               \end{pmatrix}=
                          \frac{1}{q_{1}(\zeta)}
                                    \begin{pmatrix}
                                              P_{1}(\zeta)            &  \Phi_{1}(\zeta) \\
                                    \Phi_{1}(\overline{\zeta})^{\ast} & q_{2}(\zeta)I_{r}
                                    \end{pmatrix}
                          \end{equation}
                belongs to the class $\mathcal{IP}_{k-1}^{(m+r_{1})\times(m+r_{1})}$.
                \medskip

              \textbf{(ii)} $q_{1}(\zeta)\equiv0$.
                 Then $q_{2}(\zeta)\neq0$. By Theorem \ref{th5.1},
                       \begin{equation}\label{eq7.7}
                       h_{0}(\zeta)=s_{11}-s_{12}(s_{22}-\zeta_{1}^{-1}I_{r_{1}})^{-1}s_{21}=
                               \left(
                               \begin{array}{c|c}
                               s_{11}   &   s_{12}\\
                               \hline
                               s_{21}   & s_{22}-\zeta_{1}^{-1}I_{r_{1}}
                               \end{array}
                               \right),
                       \end{equation}
                 where
                     \begin{equation}\label{eq7.8}
                        \begin{pmatrix}
                        s_{11} & s_{12} \\
                        s_{21} & s_{22}
                        \end{pmatrix}=
                            \frac{1}{q_{2}(\zeta)}
                                  \begin{pmatrix}
                                    P_{2}(\zeta)                     &  \Phi_{1}(\zeta) \\
                                   \Phi_{1}(\overline{\zeta})^{\ast} &       0_{r}
                                   \end{pmatrix}
                        \in\mathcal{IP}_{k-1}^{(m+r_{1})\times(m+r_{1})}.
                     \end{equation}
                 By the Frobenius formula \cite{uj08}, there exists a matrix
                        $$
                          \begin{pmatrix}
                           s_{22}   & I_{r_{1}} \\
                          I_{r_{1}} & \zeta_{1}I_{r_{1}}
                          \end{pmatrix}^{-1}=
                                 \begin{pmatrix}
                                 (s_{22}-\zeta_{1}^{-1})^{-1} & \star \\
                                             \star            & \star
                                 \end{pmatrix}.
                        $$
                 Then
                   \begin{multline}\label{eq7.9}
                   h_{0}(\zeta)=
                     \left(
                     \begin{array}{c|c}
                     s_{11}   &         s_{12}      \\
                     \hline
                     s_{21}   & s_{22}-\zeta_{1}^{-1}I_{r_{1}}
                     \end{array}
                     \right)=
                           \left(
                           \begin{array}{c|cc}
                           s_{11}  &   s_{12}    &      0       \\
                           \hline
                           s_{21}  &   s_{22}    &   I_{r_{1}}  \\
                             0     &   I_{r_{1}} & 0+\zeta_{1}I_{r_{1}}
                           \end{array}
                           \right).
                   \end{multline}
                 It is easy to see that the matrix of coefficients
                         \begin{equation}\label{eq7.10}
                         s_{1}(\zeta)=
                            \begin{pmatrix}
                            s_{11}(\zeta)  &  s_{12}(\zeta) &    0  \\
                            s_{21}(\zeta)  &  s_{22}(\zeta) &  I_{r_{1}}\\
                                  0        &    I_{r_{1}}   &    0
                            \end{pmatrix}
                         \end{equation}
                 of the fractional linear transformation (\ref{eq7.9})
                 is a multi-affine Cayley inner function of the class $\mathcal{IP}_{k-1}^{(m+2r_{1})\times(m+2r_{1})}$.

       In both cases \textbf{(i)}, \textbf{(ii)} we have
            \begin{equation}\label{eq7.11}
            h_{0}(\zeta)=
               \left(
               \begin{array}{c|c}
                  h_{11}(\zeta) & h_{12}(\zeta)     \\
                  \hline
                  h_{21}(\zeta) & h_{22}(\zeta)+\Lambda_{1}
               \end{array}
               \right),
            \end{equation}
       where the multi-affine coefficient matrix
                 $$
                 h_{1}(\zeta)=
                    \begin{pmatrix}
                    h_{11}(\zeta) & h_{12}(\zeta) \\
                    h_{21}(\zeta) & h_{22}(\zeta)
                    \end{pmatrix}
                 $$
       does not depend on $\zeta_{1}$ and has the form (\ref{eq7.6}) or (\ref{eq7.10}), and
           $$
           \Lambda_{1}=\zeta_{1}I_{r_{1}}\quad \mbox{or} \quad \Lambda_{1}=\mbox{diag}\{0_{r_{1}},\zeta_{1}I_{r_{1}}\}.
           $$

       \noindent
       \textbf{\emph{Step j}}, $2\leq j\leq k$.
           Suppose
               \begin{equation}\label{eq7.12}
               h_{0}(\zeta_{1},\ldots,\zeta_{k})=
                      \left(
                      \begin{array}{c|c}
                      \widehat{h}_{11}(\zeta) & \widehat{h}_{12}(\zeta)    \\
                      \hline
                      \widehat{h}_{21}(\zeta) & \widehat{h}_{22}(\zeta)+\Lambda_{j-1}
                      \end{array}
                      \right),
               \end{equation}
           where
                  \begin{equation}\label{eq7.13}
                     \begin{pmatrix}
                     \widehat{h}_{11}(\zeta) & \widehat{h}_{12}(\zeta) \\
                     \widehat{h}_{21}(\zeta) & \widehat{h}_{22}(\zeta)
                     \end{pmatrix}=
                           h_{j-1}(\zeta)
                  \end{equation}
           is a multi-affine Cayley inner function of the Pick class independent of the variables $\zeta_{1},\ldots,\zeta_{j-1}$ and
                       $$
                       \Lambda_{j-1}=\mbox{diag}\{0_{j-1},\zeta_{1}I_{r_{1}},\ldots,\zeta_{j-1}I_{r_{j-1}}\}.
                       $$
           We represent (\ref{eq7.13}) in the form
                  \begin{equation}\label{eq7.14}
                  h_{j-1}(\zeta)=\frac{\zeta_{j}P_{1}(\zeta)+P_{2}(\zeta)}{\zeta_{j}q_{1}(\zeta)+q_{2}(\zeta)}.
                  \end{equation}
           Repeating the reasoning from Step 1 for $h_{j-1}(\zeta)$, we obtain
                        \begin{equation}\label{eq7.15}
                        h_{j-1}(\zeta)=
                              \begin{pmatrix}
                              \widehat{h}_{11} & \widehat{h}_{12} \\
                              \widehat{h}_{21} & \widehat{h}_{22}
                              \end{pmatrix}=
                                   \left(
                                   \begin{array}{cc|c}
                                   a_{11}(\zeta) & a_{12}(\zeta) & a_{13}(\zeta) \\
                                   a_{21}(\zeta) & a_{22}(\zeta) & a_{23}(\zeta) \\
                                   \hline
                                   a_{31}(\zeta) & a_{32}(\zeta) & a_{33}(\zeta)+\Lambda
                                   \end{array}
                                   \right),
                        \end{equation}
           where
              $$
              \Lambda=\mbox{diag}\{0_{r_{j}},\zeta_{j}I_{r_{j}}\}.
              $$
           By Proposition \ref{pr6.1}, from (\ref{eq7.12}), (\ref{eq7.15}) we obtain
                      \begin{equation}\label{eq7.16}
                        h_{0}(\zeta_{1},\ldots,\zeta_{k})=
                              \left(
                              \begin{array}{c|cc}
                              a_{11}(\zeta)  &         a_{12}(\zeta)       & a_{13}(\zeta) \\
                              \hline
                              a_{21}(\zeta) & a_{22}(\zeta)+\Lambda_{j-1} & a_{23}(\zeta)\\
                              a_{31}(\zeta) &         a_{32}(\zeta)       & a_{33}(\zeta)+\Lambda
                              \end{array}
                              \right)
                      \end{equation}
           If necessary, then by permutation the rows and corresponding columns in the last matrix
           and introducing a new division into blocks, we get
               \begin{equation}\label{eq7.17}
               h_{0}(\zeta_{1},\ldots,\zeta_{k})=
                           \left(
                           \begin{array}{c|c}
                            h_{11}(\zeta) & h_{12}(\zeta)    \\
                            \hline
                            h_{21}(\zeta) & h_{22}(\zeta)+\Lambda_{j}
                           \end{array}
                           \right),
               \end{equation}
           where $\{h_{kl}(\zeta)\}_{k,l=1}^{2}=h_{j}(\zeta)$ is multi-affine Cayley inner function of the Pick class
           independent of the variables $\zeta_{1},\ldots,\zeta_{j}$, and
                  $$
                  \Lambda_{j}=\mbox{diag}\{0_{j},\zeta_{1}I_{r_{1}},\ldots,\zeta_{j-1}I_{r_{j-1}},\zeta_{j}I_{r_{j}}\}.
                  $$
           At the last step (for $j=k$), we obtain the matrix of coefficients of the fractional linear transformation (\ref{eq7.17}),
           which does not depend on any of the variables. Since such a Cayley inner function of the Pick class
           is a constant Hermitian matrix, we see that
                         $$
                         h_{k}(\zeta)\equiv B=
                              \begin{pmatrix}
                              B_{11} & B_{12} \\
                              B_{21} & B_{22}
                              \end{pmatrix},\quad
                                   B=B^{\ast}
                         $$
           and
             \begin{multline}\label{eq7.18}
               h_{0}(\zeta_{1},\ldots,\zeta_{k})=
                     \left(
                     \begin{array}{c|c}
                       B_{11} & B_{12}    \\
                      \hline
                       B_{21} & B_{22}+\Lambda_{k}
                     \end{array}
                     \right)=\\
                          B_{11}-B_{12}(B_{22}+\Lambda_{k})^{-1}B_{21},
             \end{multline}
           where
                $$
                \Lambda_{k}=\mbox{diag}\{0_{k},\zeta_{1}I_{r_{1}},\ldots,\zeta_{k}I_{r_{k}}\}.
                $$
           Returning to the function $f(z_{1},\ldots,z_{d})$ by identifying the variables of each $i$-th subset
           of variables with $z_{i}$, we obtain the (\ref{eq7.1}).

     \textbf{(1)} $\Rightarrow$ \textbf{(0)}.
         The trivial identity
              $$
                \begin{pmatrix}
                f(z) \\
                  0
                \end{pmatrix}=\\
                   \begin{pmatrix}
                   B_{11}  &  B_{12}  \\
                   B_{21}  & B_{22}+Z_{n}
                   \end{pmatrix}
                        \begin{pmatrix}
                             I_{m}      \\
                         -(B_{22}+Z_{n})^{-1}B_{21}
                        \end{pmatrix}
              $$
         implies
               \begin{multline}\label{eq7.19}
                f(z)=B_{11}-B_{12}(B_{22}+Z_{n})^{-1}B_{21}=(I_{m},\,-B_{21}^{\ast}(B_{22}+Z_{n})^{-1\ast})\times\\
                   \begin{pmatrix}
                   B_{11}  &  B_{12}  \\
                   B_{21}  & B_{22}+Z_{n}
                   \end{pmatrix}
                        \begin{pmatrix}
                             I_{m}      \\
                         -(B_{22}+Z_{n})^{-1}B_{21}
                        \end{pmatrix}.\qquad
               \end{multline}
         From (\ref{eq7.19}) it follows that
                      \begin{equation}\label{eq7.20}
                      \frac{f(z)-f(z)^{\ast}}{2i}=
                      B_{21}^{\ast}(B_{22}+Z_{n})^{-1\ast}\frac{Z_{n}-Z_{n}^{\ast}}{2i}(B_{22}+Z_{n})^{-1}B_{21},
                      \end{equation}
         that is, $f(z)$ is a Cayley inner function of the Pick class.
  \end{proof}

  \section{Proof of the Main Theorem}\label{s:8}

   Theorem \ref{th3.1} and \ref{th7.1} allows us to obtain a representation for a rational function $f(z)$ of the Pick class.
   \medskip

   \emph{Proof of Main Theorem}.

   \noindent
   \textbf{(0)}$\Rightarrow$\textbf{(1)}. Let be $f(z)\in\mathcal{P}_{d}^{m\times m}$. By Theorem \ref{th3.1},
               there exists $g(z_{0},z)\in\mathcal{IP}_{d+1}^{m\times m}$ such that
                       $$
                         g(i,z)=f(z).
                       $$
      By Theorem \ref{th7.1}, the function $g(z_{0},z)$ has a representation
                \begin{equation}\label{eq7.27}
                  g(z_{0},z)=B_{11}-B_{12}(B_{22}+Z)^{-1}B_{21}=
                            \left(
                            \begin{array}{c|c}
                            B_{11} & B_{12} \\
                            \hline
                            B_{21} & B_{22}+Z
                            \end{array}
                            \right).
                \end{equation}
      We rewrite the representation (\ref{eq7.27}) in a more detailed form
                    \begin{equation}\label{eq7.28}
                    g(z_{0},z)=
                            \left(
                            \begin{array}{c|ccc}
                            B_{11} & C_{12} &        C_{13}         &    C_{14}  \\
                            \hline
                            C_{21} & C_{22} &        C_{23}         &    C_{24}  \\
                            C_{31} & C_{32} & C_{33}+z_{0}I_{r_{0}} &    C_{34}  \\
                            C_{41} & C_{42} &        C_{43}         & C_{44}+\widehat{Z}_{n}
                            \end{array}
                            \right),
                    \end{equation}
      where
           $$
           \widehat{Z}_{n}=\mbox{diag}\{z_{1}I_{n_{1}},\ldots,z_{d}I_{n_{d}}\}.
           $$
      Setting in (\ref{eq7.28}) $z_{0}=i$, we obtain
                \begin{equation}\label{eq7.29}
                    f(z_{1},\ldots,z_{d})=
                            \left(
                            \begin{array}{c|ccc}
                            B_{11} & C_{12} &        C_{13}         &    C_{14}  \\
                            \hline
                            C_{21} & C_{22} &        C_{23}         &    C_{24}  \\
                            C_{31} & C_{32} & C_{33}+iI_{r_{0}} &    C_{34}  \\
                            C_{41} & C_{42} &        C_{43}         & C_{44}+\widehat{Z}_{n}
                            \end{array}
                            \right).
                \end{equation}
      The matrix $H$ of the fractional linear transformation (\ref{eq7.29}) satisfies the condition $(H-H^{\ast})/2i\geq0$, as required.
      \medskip

   \noindent
   \textbf{(1)}$\Rightarrow$\textbf{(2)}.
       If $f(z_{1},\ldots,z_{d})\in\mathcal{P}_{d}^{m\times m}$, then
               $$
               g(z_{1},\ldots,z_{d})=f(-1/z_{1},\ldots,-1/z_{d})\in\mathcal{P}_{d}^{m\times m}.
               $$
       By \textbf{(1)}, we have
                    $$
                    g(z)=
                         \left(
                         \begin{array}{c|cc}
                         H_{11}  & H_{12}  & H_{13}  \\
                         \hline
                         H_{21}  & H_{22}  & H_{23}  \\
                         H_{31}  & H_{32}  & H_{33}+\widehat{Z}_{n}
                         \end{array}
                         \right),
                    $$
       where
           $$
           \widehat{Z}_{n}=\text{diag}\{z_{1}I_{n_{1}},\ldots,z_{d}I_{n_{d}}\}.
           $$
       Then
             \begin{multline}\label{eq7.30}
              f(z_{1},\ldots,z_{d})=g(-1/z_{1},\ldots,-1/z_{d})=\\
                   \left(
                   \begin{array}{c|cc}
                   H_{11}  &              H_{12}                & H_{13}  \\
                   \hline
                   H_{21}  & (H_{22}+I_{n_{0}})-I_{n_{0}}^{-1}  & H_{23}  \\
                   H_{31}  &              H_{32}                & H_{33}-\widehat{Z}_{n}^{-1}
                   \end{array}
                   \right).
             \end{multline}
       From (\ref{eq7.30}) it follows that
                   \begin{equation}\label{eq7.31}
                   f(z)=A-B(D-Z_{n}^{-1})^{-1}C=A+BZ_{n}(I-DZ_{n})^{-1}C,
                   \end{equation}
       where
            $$
             A=H_{11},\, B=(H_{12},\,H_{13}),\, C=
                    \begin{pmatrix}
                    H_{21} \\
                    H_{31}
                    \end{pmatrix},\
                        D=
                          \begin{pmatrix}
                          H_{22}+I_{n_{0}} & H_{23} \\
                               H_{32}      & H_{33}
                          \end{pmatrix},
            $$
       and $Z_{n}=\text{diag}\{I_{n_{0}},z_{1}I_{n_{1}},\ldots,z_{d}I_{n_{d}}\}$. The inequality (\ref{eq7.22}) is obvious.
       \medskip

   \noindent
   \textbf{(2)}$\Rightarrow$\textbf{(1)}.
       If $f(z)\in\mathcal{P}_{d}^{m\times m}$, then
           $$
           g(z_{1},\ldots,z_{d})=f(-1/z_{1},\ldots,-1/z_{d})\in\mathcal{P}_{d}^{m\times m}.
           $$
       By \textbf{(2)}, we have
              $$
              g(z)=\widehat{A}+\widehat{B}Z(I-\widehat{D}Z)^{-1}\widehat{C}=
                \widehat{A}-\widehat{B}(\widehat{D}-Z^{-1})^{-1}\widehat{C},
              $$
       where
           $$
           Z=
             \begin{pmatrix}
               I_{n_{0}} &     0   \\
                  0      & \widehat{Z}_{n}
             \end{pmatrix}, \quad
                \widehat{Z}_{n}=\text{diag}\,\{z_{1}I_{n_{1}},\ldots,z_{d}I_{n_{d}}\}.
           $$
       Then
           $$
           g(z)=
              \left(
                \begin{array}{c|cc}
                H_{11} &       H_{12}          &        H_{13} \\
                \hline
                H_{21} & H_{22}-I_{n_{0}}^{-1} &        H_{23} \\
                H_{31} &       H_{32}          & H_{33}-\widehat{Z}_{n}^{-1}
                \end{array}
              \right).
           $$
       From here
               $$
               f(z)=g(-1/z_{1},\ldots,-1/z_{d})=
                      \left(
                      \begin{array}{c|cc}
                      H_{11} &       H_{12}       &        H_{13} \\
                      \hline
                      H_{21} &   \widehat{H}_{22} &        H_{23} \\
                      H_{31} &       H_{32}       & H_{33}+\widehat{Z}_{n}
                      \end{array}
                      \right),
               $$
       where $\widehat{H}_{22}=H_{22}-I_{n_{0}}^{-1}$, that is,
                  $$
                  f(z)=A-B(D+Z_{n})^{-1}C,\quad Z_{n}=\text{diag}\,\{0_{n_{0}},z_{1}I_{n_{1}},\ldots,z_{d}I_{n_{d}}\}.
                  $$

   \noindent
   \textbf{(1)}$\Rightarrow$\textbf{(3)}.
        In fact, the representation (\ref{eq7.22})
           \begin{multline}\label{eq7.33}
           f(z)=A-B(D+Z_{n})^{-1}C=\\
                   \left(
                   \begin{array}{c|cccc}
                   H_{11}    &  H_{12}   &        H_{13}          & \cdots &    H_{1,d+2}  \\
                   \hline
                   H_{21}    &  H_{22}   &        H_{23}          & \cdots &    H_{2,d+2} \\
                   H_{31}    &  H_{32}   & H_{33}+z_{1}I_{n_{1}}  & \cdots &    H_{3,d+2}  \\
                   \vdots    &  \vdots   &        \vdots          & \ddots &    \vdots    \\
                   H_{d+2,1} & H_{d+2,2} &       H_{d+2,3}        & \cdots &    H_{d+2,d+2}+z_{d}I_{n_{d}}
                   \end{array}
                   \right)
           \end{multline}
        is a representation of $f(z)$ in the form of a long-resolvent \cite{uj06,uj07}:
               $$
               f(z)=A_{11}(z)-A_{12}(z)A_{22}(z)^{-1}A_{21}(z),
               $$
        where the positive semidefinite matrix coefficients $A_{k}$ for independent variables of the matrix pencil
                   $$
                    A(z)=
                        \begin{pmatrix}
                        A_{11}(z) & A_{12}(z) \\
                        A_{21}(z) & A_{22}(z)
                        \end{pmatrix}=
                            H+z_{1}A_{1}+\cdots+z_{d}A_{d}
                   $$
        are pairwise orthogonal projectors: $A_{k}^{2}=A_{k}=A_{k}^{\ast}$, $k=1,\ldots,d$,
        $A_{j}A_{k}=A_{k}A_{j}=0$, $j\neq k$.
        \medskip

   \noindent
   \textbf{(3)}$\Rightarrow$\textbf{(0)}.
       Since $A_{22}(z)$ is an invertible matrix, we see that
               $$
               f(z)=(I_{m},\,-A_{21}^{\ast}A_{22}^{-1\ast})
                    \begin{pmatrix}
                    A_{11} & A_{12} \\
                    A_{21} & A_{22}
                    \end{pmatrix}
                        \begin{pmatrix}
                            I_{m}   \\
                        -A_{22}^{-1}A_{21}
                        \end{pmatrix}.
               $$
       From here we get
                   \begin{multline}
                   \frac{f(z)-f(z)\ast}{2i}=\\
                       (I_{m},\,-A_{21}(z)^{\ast}A_{22}(z)^{-1\ast}) \frac{H-H^{\ast}}{2i}
                         \begin{pmatrix}
                             I_{m}   \\
                         -A_{22}(z)^{-1}A_{21}(z)
                         \end{pmatrix}+\\
                              \sum_{k=1}^{d}\frac{z_{k}-\overline{z}_{k}}{2i}
                                 (I_{m},\,-A_{21}(z)^{\ast}A_{22}(z)^{-1\ast})A_{k}
                                         \begin{pmatrix}
                                                I_{m}   \\
                                         -A_{22}(z)^{-1}A_{21}(z)
                                         \end{pmatrix},
                   \end{multline}
       that is, $f(z)\in\mathcal{P}_{d}^{m\times m}$.

       The last statement of the theorem is obvious.
   \qed

 \section{Appendix. Sum-of-Squares Theorem}\label{s:Ap}

    On $\mathbb{R}^{d}$, the values of $\mathbb{C}^{m\times m}$-valued Cayley inner functions $f(z_{1},\ldots,z_{d})$
    of the Pick class are complex Hermitian matrices. For this case, the proof of Sum-of-Squares Theorem differs somewhat
    from given in \cite{uj03}, where scalar functions were considered.

    First of all, consider several auxiliary concepts and theorem.

    If $K$ is a field, then $K(x_{1},\ldots,x_{d})$ denotes the set of rational functions in variables $x_{1},\ldots,x_{d}$
    with coefficients from the fields $K$. Artin's solution to Hilbert's 17th problem for scalar polynomials is well know
    (see, for example, \cite{uj19}, Chapter XI, Corollary 3.3).
    \medskip

    \noindent
    \textbf{Artin's Theorem}.
        \emph{Let $K$ be a real fields admitting only one ordering.
        Let $f(x)\in K(x_{1},\ldots,x_{d})$ be a rational function that does not take negative values:
        $f(a)\geq 0$ for all $a=(a_{1},\ldots, a_{d})\in K^{d}$, in which $f(a)$ is defined. Then $f(x)$
        is a sum of squares in $K(x_{1},\ldots,x_{d})$.}

    \begin{remark}\label{remAp.1}
        This means: for any polynomial $g(z)\in\mathbb{R}[z_{1},\ldots,z_{d}]$ that is nonnegative
        on $\mathbb{R}^{d}$, there exists a polynomial $s(z)\in\mathbb{R}[z_{1},\ldots,z_{d}]$ such that
        $s(z)^{2}g(z)$ is an SOS polynomial.
        \medskip
        
        \noindent
        Polynomial $s(z)$ will be called \emph{an Artin denominator} of the nonnegative polynomial $g(x)\geq0,\,x\in\mathbb{R}^{d}$.
    \end{remark}

        The $\mathbb{C}^{m\times m}$-valued PSD polynomial also has Artin's denominator, which is a scalar polynomial with real coefficients.

    \begin{proposition}\label{prAp.1}
        Let $F(z)$, $z\in\mathbb{C}^{d}$ be the $\mathbb{C}^{m\times m}$-valued PSD polynomial.
        Then there exists a scalar polynomial $s(z)\in\mathbb{R}[z_{1},\ldots,z_{d}]$ such that
              \begin{equation}\label{eqAp.1}
              s(z)^{2}F(z)=M(z)M(\overline{z})^{\ast},
              \end{equation}
        where $M(z)$ is some $\mathbb{C}^{m\times k}$-valued polynomial. That is, $s(z)^{2}F(z)$
        is an SOS polynomial.
    \end{proposition}

    \begin{proof}
        Let $\Phi(x)$ be $\mathbb{R}^{m\times m}$-valued PSD polynomial. Using the Jacobi method of reducing a quadratic form to
        canonical form, Sylvester's criterion and Artin's Theorem we obtain: there exists a polynomial $s(z)\in\mathbb{R}[z_{1},\ldots,z_{d}]$
        such that $s(z)^{2}\Phi(z)=R(z)R(z)^{T}$, where $R(x)$ is some $\mathbb{R}^{m\times k}$-valued polynomial.

        Let $F(z)$ be a $\mathbb{C}^{m\times m}$-valued PSD polynomial and
            \begin{equation}\label{eqAp.2}
            A(z)=\frac{F(z)+\overline{F(\overline{z})}}{2},\quad B(z)=\frac{F(z)-\overline{F(\overline{z})}}{2i}.
            \end{equation}
        Since $F(\overline{z})=F(z)^{\ast}=\overline{F(z)}^{T}$, we see that $A(x),\,B(x)$ are $\mathbb{R}^{m\times m}$-valued polynomial
        satisfying the conditions
               \begin{equation}\label{eqAp.3}
               A(z)^{T}=A(z),\quad\text{and}\quad B(z)^{T}=-B(z).
               \end{equation}
        For every row vector $\eta=\xi_{1}+i\xi_{2}\in\mathbb{C}^{m}$, where $\xi_{1},\xi_{2}\in\mathbb{R}^{m}$, we have
                     \begin{multline}\label{eqAp.4}
                     \eta F(x)\eta^{\ast}=(\xi_{1}+i\xi_{2})(A(x)+iB(x))(\xi_{1}^{T}-i\xi_{2}^{T})=\\
                          (\xi_{1},\,\xi_{2})
                              \begin{pmatrix}
                                A(x)   & B(x)\\
                              B(x)^{T} & A(x)
                              \end{pmatrix}
                                  \begin{pmatrix}
                                  \xi_{1}^{T}\\
                                  \xi_{2}^{T}
                                  \end{pmatrix}\geq0.\quad
                    \end{multline}
        Thus,
           $$
           \Phi(x)=
              \begin{pmatrix}
                  A(x)   & B(x)\\
                B(x)^{T} & A(x)
              \end{pmatrix}
           $$
        is the $\mathbb{R}^{2m\times 2m}$-valued PSD polynomial. Therefore, there exist $s(z)\in\mathbb{R}[z_{1},\ldots,z_{d}]$
        and $\mathbb{R}^{m\times k}$-valued polynomials $R_{1}(x),\,R_{2}(x)$ for which
                 \begin{equation}\label{eqAp.5}
                 s(z)^{2}
                     \begin{pmatrix}
                       A(z)   & B(z)\\
                     B(z)^{T} & A(z)
                     \end{pmatrix}=
                          \begin{pmatrix}
                          R_{1}(z) \\
                          R_{2}(z)
                          \end{pmatrix}
                                 (R_{1}(z)^{T},\,R_{2}(z)^{T}).
                 \end{equation}
        From (\ref{eqAp.3}), (\ref{eqAp.5}) it follows that
                      \begin{equation}\label{eqAp.6}
                       R_{1}(z)R_{1}(z)^{T}=R_{2}(z)R_{2}(z)^{T}=s(z)^{2}A(z),
                      \end{equation}
                           \begin{equation}\label{eqAp.7}
                           R_{1}(z)R_{2}(z)^{T}=s(z)^{2}B(z)=-R_{2}(z)R_{1}(z)^{T}.
                           \end{equation}
        Let us $H(z)=(R_{1}(z)-iR_{2}(z))/\sqrt{2}$. Taking into account (\ref{eqAp.6}), (\ref{eqAp.7}), we obtain
                \begin{multline}\label{eqAp.8}
                \quad H(z)H(\overline{z})^{\ast}=(R_{1}(z)-iR_{2}(z))(R_{1}(z)^{T}+iR_{2}(z)^{T})/2=\\
                (R_{1}(z)R_{1}(z)^{T}+R_{2}(z)R_{2}(z)^{T})/2+i(-R_{2}(z)R_{1}(z)^{T}+R_{1}(z)R_{2}(z)^{T})/2=\\
                s(z)^{2}(A(z)+iB(z))=s(z)^{2}F(z),\quad
                \end{multline}
        i.e., $s(z)^{2}F(z)$ is the $\mathbb{C}^{m\times m}$-valued SOS polynomial.
    \end{proof}

    \begin{definition}\label{defAp.1}
          \emph{Artin's denominator $s$ of PSD not SOS polynomial $F$ is called \emph{an Artin minimal denominator},
          if a polynomial $\widehat{s}=s/s_{j}$ is not Artin's denominator of $F$ for every irreducible factor $s_{j}$ of $s$.}
    \end{definition}

    \begin{theorem}\label{thAp.1}
          Each $\mathbb{C}^{m\times m}$-valued PSD not SOS polynomial $F(z)$ has a non-constant Artin minimal denominator $s(z)$.
          The irreducible factors of Artin's minimum denominator do not change sign on $\mathbb{R}^{d}$.
    \end{theorem}

          We need some additional considerations.
          If $F(z)$ is an SOS polynomial, then $s(z)^{2}F(z)$ is also an SOS polynomial for every polynomial $s(z)\in\mathbb{R}[z_{1},\ldots,z_{d}]$.
          If PSD polynomial $F(z)$ is not representable as sum of squares of polynomials, then the question arises:
          \emph{for which $s(z)$ is the polynomial $s(z)^{2}F(z)$ also PSD not SOS?}

     \begin{proposition}\label{prAp.2} \emph{(\cite{uj15}, Lemma 2.1).}
          Let $F(x)$ be a scalar PSD not SOS polynomial and $s(x)$ an irreducible indefinite polynomial in $\mathbb{R}[x_{1},\ldots,x_{d}]$.
          Then $s^{2}F$ is also a PSD not SOS polynomial.
     \end{proposition}

     \begin{proof}
           Clearly $s^{2}F$ is PSD. If $s^{2}F=\sum_{k} h_{k}^{2}$, then for every real tuple $a$ with $s(a)=0$, it follows that $s^{2}F(a)=0$.
           This implies $h_{k}(a)^{2}=0$ $\forall k$. So on the real variety $s=0$, we have $h_{k}=0$ as well. So (see \cite{uj09}, Theorem 4.5.1)
           for each $k$, there exists $g_{k}$ so that $h_{k}=sg_{k}$. This gives $F=\sum_{k} g_{k}^{2}$, which is a contradiction.
     \end{proof}

     \begin{corollary}\label{corAp.1}
           Let $F(z)$ be a $\mathbb{C}^{m\times m}$-valued PSD not SOS polynomial, and $s(z)$ an irreducible polynomial
           in $\mathbb{R}[\,z_{1},\ldots,z_{d}\,]$. Then $s^{2}F$ is also a PSD not SOS polynomial.
     \end{corollary}

     \begin{proposition}\label{prAp.3}
           If $r(z)^{2}F(z)$ is an SOS polynomial and all irreducible factors of polynomial
           $r(z)\in\mathbb{R}[z_{1},\ldots,z_{d}]$ are indefinite, then $F(z)$ is also an SOS polynomial.
     \end{proposition}

     \begin{proof}
           Suppose $F(z)$ is a PSD not SOS polynomial. Let
               $$
               r(z)=r_{1}(z)\cdots r_{k}(z)
               $$
           be the decomposition of $r(z)$ into irreducible factors. Successively applying Corollary \ref{corAp.1} to the polynomials
                       $$
                       F_{1}=r_{1}^{2}F(z),\,F_{2}=r_{2}^{2}F_{1},\,\ldots\,,F_{k}=r_{k}^{2}F_{k-1}
                       $$
           we get $F_{k}=r^{2}F$ is PSD not SOS polynomial. Contradiction.
      \end{proof}

      \noindent
      \emph{Proof of Theorem} \ref{thAp.1}.
           by Proposition \ref{prAp.1}, for $\mathbb{C}^{m\times m}$-valued PSD polynomial $F(z)$  there exists
           Artin's denominator $a(z)\in\mathbb{R}[z_{1},\ldots,z_{d}]$ for which $a(z)^{2}F(z)$ is an SOS polynomial.
           Each irreducible factor of $a(z)$ is either indefinite or does not change sign on $\mathbb{R}^{d}$.
           Then $a(z)=r(z)s(z)$, where all irreducible factors of $s(z)$ do not change sign on $\mathbb{R}^{d}$,
           and the irreducible factors of $r(z)$ are indefinite. Let us $F_{1}(z)=s(z)^{2}F(z)$.
           By condition, $r(z)^{2}F_{1}(z)=a(z)^{2}F(z)$ is the SOS polynomial. Since all irreducible factors of  $r(z)$
           are indefinite, we see that $F_{1}(z)$ is an SOS polynomial (Proposition \ref{prAp.3}).
           Then $s(z)$ is also Artin's denominator for $F(z)$.
           Let $s_{j}(z)$ be some irreducible factor of $s(z)$. If $s(z)/s_{j}(z)$ remains the Artin denominator of $F$,
           then the factor $s_{j}(z)$ is removed from $s(z)$. Removing all ``excess" irreducible factors from the polynomial
           $s(z)$, we obtain Artin's denominator with the required properties.
           \qed
           \medskip

           \noindent
           To prove Sum-of-Squares Theorem, we need some statements from \cite{uj03}:

      \begin{proposition}\label{prAp.4} \emph{(\cite{uj03}, Proposition 3.4).}
           Let $s(z)\in\mathbb{R}[z_{1},\ldots,z_{z}]$ be a irreducible polynomial that does not change sign on $\mathbb{R}^{d}$.
           If $\partial s(z)/\partial z_{1}\neq0$, then there exists a point $z'=(z'_{1},x'_{2},\ldots,x'_{d})$, $\emph{Im}\,z'_{1}>0$,
           $x'_{2},\ldots,x'_{d}\in\mathbb{R}$ such that $s(z')=0$.
           Moreover, there exists a point  $z''\in\Pi^{d}$ for which $s(z'')=0$.
      \end{proposition}

        Let $\mathcal{Z}(h)=\{z\in\mathbb{C}\mid h(z)=0\}$ be a zero set of the polynomial $h$.

      \begin{proposition}\label{prAp.4} \emph{(\cite{uj03}, Proposition 7.1).}
          Let $s(z), h(z)\in\mathbb{R}[z_{1},\ldots,z_{d}]$ be coprime polynomials and $s(z'_{1},z')=h(z'_{1},z')=0$ for fixed
          $z'_{1}\in\mathbb{C}$, $z'\in\mathbb{K}^{d-1}$, where $\mathbb{K}=\mathbb{R}$ or $\mathbb{C}$.
          Put $\Omega=\Omega_{1}\times \Omega_{\mathbb{K}}$, where $\Omega_{1}\subset\mathbb{C}$ is some neighborhood of  $z'_{1}$,
          and $\Omega_{\mathbb{K}}\subset\mathbb{K}^{d-1}$ is a some neighborhood of $z'$.
          Then neither of the set $\mathcal{Z}(s)\cap\Omega$ and $\mathcal{Z}(h)\cap\Omega$ is subset of the other.
      \end{proposition}

      \begin{remark}\label{reAp.2}
          This means: if $s(z'_{1},z')=h(z'_{1},z')=0$, then in any neighborhood $\Omega=\Omega_{1}\times \Omega_{\mathbb{K}}$
          of the point $(z'_{1},z')$ there is a point $(z''_{1},z'')$ for which $s(z''_{1},z'')=0$ and $h(z''_{1},z'')\neq0$.
      \end{remark}

      \begin{theorem}\label{thAp.2} \emph{(\cite{uj03}, Theorem 6.1).}
           Let $f(z)=p(z)/q(z)$ be rational function with real coefficients and $s(z)\in\mathbb{R}[z_{1},\ldots,z_{d}]$. If
              $$
              s(z)^{2}W_{z_{1}}[q,p]=s(z)^{2}\left(q(z)\frac{\partial p(z)}{\partial z_{1}}-p(z)\frac{\partial q(z)}{\partial z_{1}}\right)=H(z)H(z)^{T}
              $$
           is the SOS polynomial, then there exist a real symmetric matrices $A_{j}$, $j=0,1\ldots,d$,
           where $A_{1}$ is positive semidefinite such that
                 \begin{equation}\label{eqAp.9}
                 f(z)=\frac{\Psi(\zeta)}{q(\zeta)s(\zeta)}(A_{0}+z_{1}A_{1}+\cdots+z_{d}A_{d})\frac{\Psi(z)^{T}}{q(z)s(z)},\quad\zeta,z\in\mathbb{C}^{d},
                 \end{equation}
                      \begin{equation}\label{eqAp.10}
                      W_{z_{1}}[q,p]=\frac{\Psi(z)}{s(z)}A_{1}\frac{\Psi(z)^{T}}{z(z)}=\frac{H(z)}{s(z)}\frac{H(z)^{T}}{s(z)},
                      \end{equation}
           where $\Psi(z)= (z^{\alpha_{1}},\ldots,z^{\alpha_{N}})$ is the row vector of all monomials
           $z^{\alpha_{i}}=z_{1}^{\delta_{1}}\cdots z_{d}^{\delta_{d}}$ satisfies the conditions
                $$
                \emph{deg}\,z^{\alpha_{i}}\leq \emph{deg}\,f+\emph{deg}\,s(z),\quad
                \emph{deg}_{z_{1}}z^{\alpha_{i}}\leq \emph{deg}_{z_{1}}f+\emph{deg}_{z_{1}}s(z).
                $$
      \end{theorem}

   \emph{Proof of Sum-of-Squares Theorem}.
   \smallskip

   Let $f(z)=P(z)/q(z)$ be a $\mathbb{C}^{m\times m}$-valued rational Cayley inner function of the Pick class.
   Without loss of generality, we can assume that $\overline{q(\overline{z})}=q(z)$ and $P(\overline{z})^{\ast}=P(z)$.
   By Proposition \ref{pr3.1}, all Wronskians $W_{z_{k}}[q,P]$, $k=1,\ldots,d$ are $\mathbb{C}^{m\times m}$-valued PSD polynomials.
   
   Suppose one of the Wronskians, say $W_{z_{1}}[q,P]$, is a PSD not SOS polynomial. By Proposition \ref{prAp.1},
   there exists Artin's denominator $s(z)\in\mathbb{R}[z_{1},\ldots,z_{d}]$ such that
        \begin{equation}\label{eqAp.11}
        s(z)^{2}W_{z_{1}}[q,P]=M(z)M(\overline{z})^{\ast},
        \end{equation}
   where $M(z)$ is some $\mathbb{C}^{m\times k}$-valued polynomial. By Theorem \ref{thAp.1}, we can assume that $s(z)$ is Artin's
   minimal denominator. 
   Let $s_{0}(z)\in\mathbb{R}[z_{1},\ldots,z_{d}]$ be an irreducible factor of $s(z)$ different from the constant.
   Among the elements $m_{kl}(z)$ of the matrix $M(z)$ there is a polynomial
            \begin{equation}\label{eqAp.12}
            m_{ij}(z)=a_{ij}(z)+ib_{ij}(z),
            \end{equation}
   where $a_{ij}(z),b_{ij}(z)\in\mathbb{R}[z_{1},\ldots,z_{d}]$, for which $s_{0}(z)$ is not a divisor.
     Indeed, if $s_{0}(z)$ is a divisor of all elements of $m_{kl}(z)=a_{kl}(z)+ib_{kl}(z)$, then $s_{0}(z)$ is the divisor
   of $a_{kl}(z)$ and $b_{kl}(z)$. Therefore, it is the divisor of $\overline{m_{kl}(\overline{z})}=a_{kl}(z)-ib_{kl}(z)$.
   Then $s(z)/s_{0}(z)$ is also the Artin denominator for $W_{z_{1}}[q,P]$, which contradicts the minimality of $s(z)$.
   
    The element $m_{ij}(z)$ (\ref{eqAp.12}) belongs to the $i$-th row $M_{i}(z)$ of the matrix $M(z)$. We represent it in the form
         $$
         M_{i}(z)=A(z)+iB(z),
         $$
    where the elements of the row vector $A(z),B(z)$ are polynomials with real coefficients.
    
    Note that the $i$-th diagonal element $f_{ii}(z)=p_{ii}(z)/q(z)$ of the original function $f(z)$ is a scalar rational Cayley inner
    function of the Pick class with real coefficients. From (\ref{eqAp.11}) we get the polynomial
            $$
            s(z)^{2}W_{z_{1}}[q,p_{ii}]=M_{i}(z)M_{i}(\overline{z})^{\ast}
            $$
    with real coefficients. From here
        \begin{equation}\label{eqAp.13}
        s(z)^{2}W_{z_{1}}[q,p_{ii}]=(A(z)+iB(z))(A(z)-iB(z))^{T}=H(z)H(z)^{T},
        \end{equation}
    where $H(z)=(A(z),\,B(z))=(h_{1}(z),\ldots,h_{2k}(z))$. The irreducible factor $s_{0}(z)$ is not a divisor of at least one of
    the elements $a_{ij}(z)$, $b_{ij}(z)$ (\ref{eqAp.12}) of $H(z)$. Without loss of generality, we can assume that
    $a_{ij}(z)=h_{1}(z)$ and $s_{0}(z)$ are coprime polynomial.
    
    By Theorem \ref{eqAp.1}, for $f_{ii}(z)=p_{ii}(z)/q(z)$ there exist a real symmetric matrices $A_{j}$, $j=0,1\ldots,d$,
    where $A_{1}$ is positive semidefinite, such that
          \begin{equation}\label{eqAp.14}
          f_{ii}(z)=\frac{\Psi(\zeta)}{q(\zeta)s(\zeta)}(A_{0}+z_{1}A_{1}+\cdots+z_{d}A_{d})\frac{\Psi(z)^{T}}{q(z)s(z)},
          \quad\zeta,z\in\mathbb{C}^{d},
          \end{equation}
                \begin{equation}\label{eqAp.15}
                W_{z_{1}}[q,p_{ii}]=\frac{\Psi(z)}{s(z)}A_{1}\frac{\Psi(z)^{T}}{z(z)}=\frac{H(z)}{s(z)}\frac{H(z)^{T}}{s(z)}=
                \sum_{j=1}^{2k}\frac{h_{j}(z)^{2}}{s(z)^{2}}.
                \end{equation}
    Setting $\zeta=\overline{z}$, from (\ref{eqAp.14}), (\ref{eqAp.15}) we obtain
          \begin{equation}\label{eqAp.16}
          \mbox{Im}\,f_{ii}(z)=\mbox{Im}\,z_{1}\sum_{j=1}^{2k}\frac{|h_{j}(z)|^{2}}{|q(z)s(z)|^{2}}+
          \sum_{k=1}^{d}\mbox{Im}\,z_{k}\frac{\Psi(z)A_{k}\Psi(z)^{\ast}}{|s(z)q(z)|^{2}}.
          \end{equation}
          
    \noindent
    Note that the expression on the left in (\ref{eqAp.15}) is a polynomial. Therefore, $\sum_{j=1}^{2k}h_{j}(z)^{2}/s(z)^{2}$
    is the polynomial. In contrast to (\ref{eqAp.15}), the expression $\sum_{j=1}^{2k}|h_{j}(z)|^{2}/|s(z)|^{2}$ in (\ref{eqAp.16})
    may not be a polynomial, which is important in our case.
      
    For the Artin denominator $s(z)$ there are 2 possibilities: \textbf{(a)} there exists a irreducible factor $s_{0}(z)$
    of the polynomial $s(z)$ such that $\partial s_{0}/\partial z_{1}\neq0$, \textbf{(b)} $\partial s(z)/\partial z_{1}\equiv0$.
    \medskip  

    \emph{Case} \textbf{(a)}. $\partial s_{0}/\partial z_{1}\neq0$.
     
    \noindent
    The polynomial  $s_{0}(z)$ depends on at least 2 variables.
    By Proposition \ref{prAp.3}, there exists a point $z'=(z'_{1},x'_{2},\ldots,x'_{d})$, $\mbox{Im}\,z'_{1}>0$,
    $x'_{2},\ldots,x'_{d}\in\mathbb{R}$ such that $s(z')=s_{0}(z')=0$.
    Without loss of generality, we can assume that $h_{1}(z')\neq0$. Indeed, $h_{1}(z),\,s_{0}(z)$ are coprime polynomials.
    If $h_{1}(z')=s_{0}(z')=0$, then by Proposition \ref{prAp.4} there is a point $\widehat{z}'=(\widehat{z}'_{1},\widehat{x}')$,
    $\mbox{Im}\,\widehat{z}'_{1}>0$, $\widehat{x}'\in\mathbb{R}^{d-1}$ for which $h_{1}(\widehat{z}')\neq0$, $s_{0}(\widehat{z}')=0$.
    
    Since $q(z)$ is a real stable polynomial, we see that for fixed $x'\in\mathbb{R}^{d-1}$ the equation $q(z_{1},x')=0$ has only real roots.
    Therefore, $q(z'_{1},x')\neq0$ and rational function $f_{ii}(\zeta,x')=p_{ii}(\zeta,x')/q(\zeta,x')$ is holomorphic at $\mbox{Im}\,\zeta>0$.
    
    On the other hand, since $s(z'_{1},x')=0$ and $h_{1}(z'_{1},x')\neq0$, then from (\ref{eqAp.16}) we obtain
         $$
         \lim_{\zeta\rightarrow z'_{1}}\mbox{Im}\,f_{ii}(\zeta,x')=
         \lim_{\zeta\rightarrow z'_{1}}\mbox{Im}\,\zeta\cdot\sum_{j=1}^{2k}\frac{|h_{j}(\zeta,x')|^{2}}{|s(\zeta,x')q(\zeta,x')|^{2}}=+\infty,
         $$
    that is, $f_{ii}(\zeta,x')$ has a pole at $\zeta=z_{1},\,\mbox{Im}\,z_{1}>0$. A contradiction.
    \medskip
    
    \emph{Case} \textbf{(b)}.
    
    \noindent
    If $\partial s/\partial z_{1}\equiv0$, then for some $k\neq1$ there exists a irreducible factor $s_{0}(z)$, $z\in\Pi^{d-1}$
    such that $\partial s_{0}/\partial z_{k}\neq0$ holds. In addition, $s_{0}(z)$ depends on at least 2 variables.
    By proposition \ref{prAp.3}, there exists a point $z''\in\Pi^{d-1}$ for which $s_{0}(z'')=0$.
    Without loss of generality, for fixed $x_{1}\in\mathbb{R}$, we can assume that $h_{1}(x_{1},z'')\neq0$. 
    Indeed, if $h_{1}(x_{1},z'')=s_{0}(z'')=0$, then by Proposition \ref{prAp.4} there is a point $(x''_{1},\widehat{z}'')$.
    $x''_{1}\in\mathbb{R}$, $\widehat{z}''\in \Pi^{d-1}$, for which $h_{1}(x''_{1},\widehat{z}'')\neq0$, $s_{0}(\widehat{z}'')=0$.
    
    Since $q(z)$ is a real stable polynomial, we see that $q(x_{1},z'')\neq0$ for $z''\in\Pi^{d-1}$.
    Then there exists neighborhood $\Omega''\subset\Pi^{d-1}$ of point $z''\in\Pi^{d-1}$ such that $f_{ii}(x_{1},z)$ is holomorphic 
    for $z\in\Omega''$ and hence bounded. Then for some $0<C<+\infty$, using (\ref{eqAp.16}), we have
            \begin{equation}\label{eqAp.17}
            \left|\mbox{Im}\,f_{ii}(x_{1},z)\right|=
               \left|\sum_{k=2}^{d}\mbox{Im}\,z_{k}\frac{\Psi(x_{1},z)A_{k}\Psi(x_{1},z)^{\ast}}{|q(x_{1},z)s(z)|^{2}}\right|\leq C,\quad z\in\Omega''.
            \end{equation}
    Since $q(x_{1},z)\neq0$, we see that $q(z_{1},z)\neq0$ holds for $z\in\Omega''\subset\Pi^{d-1}$ and all $z_{1}$ from some complex neighborhood
    $\Omega_{1}\subset\mathbb{C}$ of the point $x_{1}\in\mathbb{R}$. Then 
         $$
         |q(z'_{1},z)s(z)|\geq\varepsilon>0,\quad \mbox{for}\, z\in\Omega''\enskip \mbox{and fixed}\enskip z'_{1}\in\Omega_{1},\, \mbox{Im}\,z'_{1}>0.
         $$
    Thus, there exists $0<C_{1}<+\infty$ such that
               \begin{equation}\label{eqAp.18}
               \left|\sum_{k=2}^{d}\mbox{Im}\,z_{k}\frac{\Psi(z'_{1},z)A_{k}\Psi(z'_{1},z)^{\ast}}{|q(z'_{1},z)s(z)|^{2}}\right|\leq C_{1},
               \enskip z\in\Omega'',\enskip\mbox{Im}\,z'_{1}>0.
               \end{equation}
    On the other hand, for a point $(z'_{1},z)\in\Omega_{1}\times\Omega''$, from (\ref{eqAp.16}) we obtain
          \begin{multline}\label{eqAp.19}
          \mbox{Im}\,f_{ii}(z'_{1},z)=\mbox{Im}\,z'_{1}\sum_{j=1}^{2k}\frac{|h(z'_{1},z)|^{2}}{|q(z'_{1},z)s(z)|^{2}}+\\
          \sum_{k=2}^{d}\mbox{Im}\,z_{k}\frac{\Psi(z'_{1},z)A_{k}\Psi(z'_{1},z)^{\ast}}{|q(z'_{1},z)s(z)|^{2}},\,\mbox{Im}\,z'_{1}>0,\,z\in\Pi^{d-1}.
          \end{multline}
    According to (\ref{eqAp.18}), the second term in (\ref{eqAp.19}) is bounded for $z\in\Omega''$.
    Since $h_{1}(z'_{1},z'')\neq0$ and $s(z'')=s_{0}(z'')=0$, we see that the first term increases indefinitely at $z\rightarrow z''\in\Pi^{d-1}$,
    which contradict holomorphy of $f(z)$ in $\Pi^{d}$.
    
    We assumed that the PSD polynomial $W_{z_{1}}[q,P]$ is not SOS, and we obtained a contradiction.
    Then $W_{z_{1}}[q,P]$ is a SOS polynomial.
    \qed

\end{document}